\documentclass[11pt, twoside, a4paper]{article}
\usepackage[english]{babel}
\usepackage{amsmath,amssymb,amsthm,epsfig}
\newtheorem{theorem}{Theorem}
\newtheorem{lemma}[theorem]{Lemma}
\newtheorem{proposition}[theorem]{Proposition}
\newtheorem{corollary}[theorem]{Corollary}

\def\supp{\mathop{\rm supp}\nolimits}

\renewenvironment{proof}[1][.]{%
\bigskip\noindent{\bf Proof#1 }}{%
\hfill$\blacksquare$\bigskip}
\title{Generic properties of Lagrangians\\ on surfaces: the
Kupka-Smale theorem}
\author{Elismar R. Oliveira\footnote{Supported by CAPES,  scholarship.}}
\date{\today}
\begin{document}
%#######################################################################################
\maketitle

\centerline{\scshape \footnotesize
Elismar R. Oliveira}
%\medskip
{\footnotesize \centerline{ Departmento de Matem\'atica}
\centerline{ Universidade Federal do Rio Grande do Sul }
\centerline{ Av. Bento Goncalves, 9500,  91509-900 }
\centerline{Porto Alegre, RS - Brasil} }\medskip

\begin{abstract}
{\footnotesize We consider generic properties of
Lagrangians. Our main result is the Theorem of Kupka-Smale, in the
Lagrangian setting, claiming that, for a convex and superlinear 
Lagrangian defined in a compact surface, for each  $k\in
\mathbb{R}$, generically, in Ma\~n\'e's sense, the energy 
level, $k$, is regular and all periodic orbits, in this level, 
are nondegenerate at all orders, that is, the linearized 
Poincar\'e map, restricted  to this energy
level, does not have roots of the unity as eigenvalues. Moreover, 
all heteroclinic intersections in this level are transversal. All
the results that we present here are true in dimension $n \geq 2$, 
except one (Theorem~\ref{PerturbacaoLocaldaOrbita}), whose proof we are able
to obtain just for dimension 2.}
\end{abstract}
%########################################################################################
%########################################################################################

%########################################################################################
\section{Introduction}
%########################################################################################
Our main purpose here is to obtain generic properties, in the sense
of Ma\~n\'e ~(see~\cite{GoReMin},~\cite{ManeGenProp}), for a convex
and superlinear Lagrangian, in a fixed smooth and compact surface
$M$ without boundary.

Our main result is the  Kupka-Smale Theorem ,
claiming that, for each, $k\in \mathbb{R}$, 
generically, this level is regular, all periodic orbits in this level are
nondegenerate at all orders and all heteroclinic intersections
in this level are transversal. Where, Nondegeneracy  of order $m$, 
means that,  the linearized Poincar\'e map, restricted to this energy level, 
does not have $m$-roots of the unity as eigenvalues.

In the proof  of the  Kupka-Smale Theorem , we will use the
Nondegeneracy Lemma (Lemma~\ref{NondegeneracyLemma}) and a  Perturbation Lemma (Lemma~\ref{PertTransvVariedEstInst}) for 
Lagrangian submanifolds in order to
get the transversality of heteroclinic 
intersections.

The  Kupka-Smale Theorem , in this formulation, resembles the Bumpy
Metrics Theorem, for geodesic flows, formulated by R. Abraham in
1968, and proved by D. V. Anosov in  1983
(Anosov,~\cite{Anos}).

The work of W. Klingenberg and F. Takens~\cite{KgTk} in the Bumpy
Metrics Theorem proof was corrected by Anosov~\cite{Anos} using an
induction method similar to the one used by M. Peixoto~\cite{Peix} in
the proof of the classical Kupka-Smale Theorem.

In this work, we will employ the same techniques used by  Anosov~\cite{Anos},
adapted to perturbations by potentials. In order to apply to the case of perturbation by
potentials, it is necessary to introduce a 
modification of the standard argument in the Control Theory for
differential equations, initially used by Klingenberg~\cite{Kg}, for
geodesic flow perturbation setting,  by J. A.
Miranda~\cite{Jan} for magnetic flows on surfaces and by
Contreras~\cite{GonDomtd}  for the proof of Franks' Lemma for
geodesic flows.  In this case we do not  have especial coordinates, like Fermi
coordinates,  as in ~\cite{Kg},~\cite{Jan} and ~\cite{GonDomtd}, thus
we introduce a new method without the use of tubular neighborhoods,
that solves this trouble. In the beginning of the proof
we use an argument similar to the one used by Robinson~\cite{RobI}, Lemma
19, Pg. 592, but the proof is quite different.

Observe that the generic properties in Ma\~n\'e's sense cannot be
obtained  from the pioneering work  of Robinson in the general Hamiltonian setting
(see~\cite{RobI} and~\cite{RobII}) because the set of all
Hamiltonians is bigger than the set of all potentials in $M$.

The transversality is the easy part. Here we follow the approach of
Contreras \& Paternain~\cite{GonzaloPaternGenGeodFlowPositEntropy},
Lemma 2.6 or J.A. Miranda~\cite{Jan}, Lemma 3.9. The problem in this
case is to construct an explicit potential that represents the
perturbation.

The main obstacle in the proof of the Kupka-Smale Theorem  in dimension $n >2$ is
the nature of the perturbation constructed. As in
Contreras~\cite{GonDomtd}, Lemma 7.3 and 7.4,  we need to solve a
matrix equation in the Lie algebra of the symplectic group $S_p(n)=\{$Symplectic
matrices $2n \times 2n\}$. The solubility of this equation is
strongly related with the existence of repeated eigenvalues 
of the matrix $H_{pp}$ in local coordinates. The problem
is that we can not change this characteristic  by adding a
potential. Moreover, the equation involved is very complex too. 

However, we point out that the Kupka-Smale Theorem  in dimension 2 is a strong result in the study of generic Lagrangians because it works below the critical level. More than that, it can be combined with other results on the structure of Aubry-Mather sets in surfaces, like Haeflinger theorem for Mather measures with rational homology in an orientable surface, claiming that such measures are supported in periodic orbits, and  results on the nonexistence of conjugated points in the Aubry-Mather sets from Contreras and Iturriaga~\cite{GoReConvex}, in order to guarantee  the hyperbolicity of the periodic orbits in this set.

%###################################################################################
\section{The  Kupka-Smale Theorem }
\vspace{0.3cm}
%###################################################################################

We consider $(M;g)$ a, $n$- dimensional, smooth and compact, Riemannian
manifold without boundary, $L: TM \rightarrow
\mathbb{R} $, a Lagrangian in $M$, convex and fiberwise
superlinear (see~\cite{GoReMin} to definitions) and $H: T^{*} M
\rightarrow \mathbb{R}$ the associated Hamiltonian obtained by
Legendre transform.

In the study of generic properties of Lagrangians we use the
concept of genericity due to Ma\~n\'e. The idea is that, 
the properties studied in the Aubry-Mather theory become much more
strongest in this generic setting. For more details and applications
see \cite{GoReConvex}, \cite{GoReMin}, \cite{Jan} and
\cite{ManeGenProp}.

We will say that a property $\mathcal{P}$ is generic, in Ma\~n\'e's
sense,  for $L$,  if there exists a generic set
$\mathcal{O} \subset C^{\infty}(M;\mathbb{R})$, in $C^{\infty}$ topology,
such that, for all $f \in \mathcal{O}$, $L+f$ has the property $\mathcal{P}$.

Consider $E_{L}(x,v)=\frac{\partial L}{\partial v} (x,v)\cdot v -
L(x,v)$ the energy function associated to $L$ and $\varepsilon
^{k}_{L}=\{(x,v) \in TM
\mid E_{L}(x,v)=k\}$ the set of all points in the energy level $k$.

Let $\theta \in TM$ be a periodic
point of positive period, $T_{min}$ of the Euler-Lagrange 
flow $ \phi_ {t} ^ {L}: TM \to  TM$. 
Fixed a local section transversal to this flow, $\Sigma$
contained in the energy level of $\theta$, there
exists a smooth function $\tau : U \subset \Sigma \to
\mathbb{R}$, such that, $\tau(\theta)=T_{min} $ which is the
time of first return to $\Sigma$,
such that the map $P(\Sigma, \theta): U \to
\Sigma$ given by
$$P(\Sigma, \theta)(\theta)=\phi_{\tau(\theta)}^{L}(\theta)$$
is a local diffeomorphism and $\theta$ is a fix point of $P(\Sigma, \theta)$.
This map is called Poincar\'e  first return map.
We will say that $ \theta $
(or the orbit of $ \theta $) is a
nondegenerate orbit of order $m \geq 1$ for $L$ if
$$Ker((d_{\theta} P(\Sigma, \theta))^{m}-Id)=0.$$
The property of Nondegeneracy of order $m$ means that $d_{\theta} P(\Sigma, \theta)$
does not have $m$-roots of the unity as eigenvalues.

If we are interested in the Hamiltonian viewpoint of the described
Lagrangian dynamics, then we consider the Hamiltonian $H$ associated to $L$  by the
Legendre transform in the speed, that is, 
$$\displaystyle H(x, p) = \max_ {v \in T_ {x} M} \{pv - L (x, v) \}.$$
Let $X^ {H} $ be the Hamiltonian field, which is the unique field
$X^ {H} $ in $T^{*} M$ such that $ \omega_ {\vartheta} (X^{H}
(\vartheta), \xi) = d_ {\vartheta} H \xi$ for all $ \xi \in T_
{\vartheta} T^{*} M$. In the local coordinates $ (x, p) $, 
$X^ {H} = H_{p} \, \frac{\partial \,}{\partial x} -
H_{x} \, \frac{\partial \,}{\partial p} $.  We denote by
$ \psi_ {t} ^ {H}: T^ {*} M \to   T^{*} M$ the
flow in $T^{*}M$ associated with the Hamiltonian field $X^ {H}: T^
{*} M \rightarrow TT^ {*} M$. This flow preserves the canonical symplectic
form $\omega$.  Since $L$ is a convex and superlinear Lagrangian  we have that $H$ is a
convex and superlinear Hamiltonian. Using the Legendre transform
\begin{center} $p=L_v (x, v) $ and $v=H_p (x, p) $, \end{center}
we have that, $H_ {pp} (x, p) $ is positive defined in $T^ {*}_{x} M$,
uniformly in $x \in M$. Observe that the Legendre transform associates the energy
level $\varepsilon_{k}^{L} $ with the level set $H^{- 1}(k)$ of $H$.
From the conjugation property between Lagrangian and Hamiltonian viewpoint, the nondegeneracy of an orbit is equivalent in both senses.

One can prove that the restriction of the symplectic form $\omega$
to $T_ {\vartheta} \Sigma$ is nondegenerate and closed form,
therefore the Poincar\'e map is symplectic.Moreover, 
$$d_ {\vartheta} \psi_ {T_ {min}} ^ {H} (\xi) = - d_
{\vartheta} \tau (\xi) X^ {H} + d_ {\vartheta} P (\Sigma, \vartheta)
(\xi), \; \forall \xi \in T_ {\vartheta} \Sigma.$$ Therefore we have
that \begin{center} \( d_ {\vartheta} \psi_ {T_ {min}} ^ {H} \mid_
{T_ {\vartheta} H^ {- 1} (k)} = \left [ \begin{matrix} 1 & d_
{\vartheta} \tau   \\ 0   & d_ {\vartheta} P (\Sigma, \vartheta)
\end{matrix} \right], \) \end{center} in general for $T=m T_
{min} $ \\ $d_ {\vartheta} \psi_ {T}^{H} (\xi) = -d_ {\vartheta}
\tau (\sum_ {i=0} ^ {M-1} d_ {\vartheta} P (\Sigma, \vartheta)^ {i})
(\xi) X^{H} + d_ {\vartheta} P (\Sigma, \vartheta)^{m} (\xi), \;
\forall \xi \in T_ {\vartheta} \Sigma$.

So, the condition of that $ \vartheta $ is nondegenerate of
order $m \geq 1$ is equivalent to say that the algebraic
multiplicity of $ \lambda=1$ as eigenvalue of $d_ {\vartheta} \psi_
{T}^{H} \mid_{T_ {\vartheta} H^ {- 1} (k)} $ is equal to 1, because
the characteristic polynomials are related by
$\mathfrak {p}_{d_{\vartheta} \psi_{T}^{H}} (\lambda) = (1- \lambda)
\cdot \mathfrak{p}_{d_{\vartheta} P(\Sigma, \vartheta)^{m}}
(\lambda) $.

Our main result is the Kupka-Smale Theorem  that relies the 
Bumpy Metrics Theorem proved by Anosov~\cite{Anos}, but here
for the Lagrangian setting.

%###################################################################################

We state our main result just when $dim(M)=2$  because we do not know
the proof for Theorem~\ref{PerturbacaoLocaldaOrbita} when $dim (M)=n
\geq 3$. This problem remains an open question.

\begin{theorem}\label{KS}\textnormal{(Kupka-Smale Theorem)}
 Suppose $dim (M)=2$. Let $L: TM \to
\mathbb{R} $, be a Lagrangian in $M$, convex and fiberwise
superlinear. Then, for each  $k \in \mathbb{R}$, the property\\
 i) $\varepsilon ^{k}_{L}$ is regular;\\
 ii) Any periodic orbit in the level $\varepsilon ^{k}_{L}$ is
 nondegenerate for all orders;\\
 iii) All heteroclinic intersections, in this level, are transversal.\\
 is generic for $L$.
\end{theorem}
%###################################################################################
\section{Proofs of the main results}

Given $k \in \mathbb{R}$, we define the set of the regular potentials for $k$, as being
$$ \mathcal{R}(k)=\{ \;f \in C^{\infty}(M;\mathbb{R}) \mid
\varepsilon ^{k}_{f}:=(H+f)^{-1}(k) \; is \; regular \; \},$$
 where $H$ is the associated Hamiltonian.

\begin{lemma}\label{compacidade}
Consider  $k \in \mathbb{R}$ and $f_{0} \in C^{\infty}(M;\mathbb{R})$.
For each sequence
$f_{n} \to f_{0}$ in $C^{\infty}(M;\mathbb{R})$ topology and points
$\vartheta_{n}=(x_{n},p_{n}) \in \varepsilon ^{k}_{f_{n}}$
there exists a subsequence $\vartheta_{n_{i}} \to \vartheta_{0} \in \varepsilon ^{k}_{f_{0}}$.
\end{lemma}

In fact, this lemma is an easy consequence of the compactness of the energy level.

\begin{theorem}\label{PotRegAbertoDenso}\textnormal{(Regularity of the energy level)}
 Given $k \in \mathbb{R}$, the subset
 $\mathcal{R}(k)$ is open and dense in
 $C^{\infty}(M;\mathbb{R})$.
\end{theorem}
\begin{proof}
The openness of the set $\mathcal{R}(k)$ follows directly of the
Lemma~\ref{compacidade}. In order to obtain the density of $\mathcal{R}(k)$ in
$C^{\infty}(M;\mathbb{R})$, consider $f_{0} \in C^{\infty}(M;\mathbb{R})$
and $\mathcal{U}$, a neighborhood  that contains a ball of radius $\varepsilon > 0$, 
and center, $f_{0}$.  We claim that $\mathcal{U} \cap \mathcal{R}(k) \neq \varnothing$. In fact, if it is not the case, we can achieve a contradiction by considering the Hamiltonian $H_{\delta}:=H + (f_{0}+ \delta)$, with $\delta \in (0,\varepsilon)$.
\end{proof}

%###################################################################################
{\large \textbf{The  Nondegeneracy Lemma }}\\
%###################################################################################

Given $k \in \mathbb{R}$ and $0< a \leq b \in \mathbb{R}$,
we define the set $\mathcal{G}_{k}^{a,b} \subseteq \mathcal{R}(k)$
as   $\mathcal{G}_{k}^{a,b}= \{f \in \mathcal{R}(k) \mid $
all periodic points  $\vartheta \in (H+f)^{-1}(k)$, 
with $T_{min}(\vartheta) \leq a$ are nondegenerate of order $m$ 
for $H+f$, \, $\forall \, m \leq \frac{b}{T_{min}} \}$.  Observe that, 
   if we have $a,a',b,b' \in \mathbb{R}$, such that,
   $0 < a, a'  < \infty$, $a \leq a'$ and $ b \leq b'$, then
   $\mathcal{G}_{k}^{a',b'} \subseteq \mathcal{G}_{k}^{a,b}$.
%###################################################################################
We define $\displaystyle
\mathcal{G}(k)=\bigcap_{n=1}^{+ \infty} \mathcal{G}_{k}^{n,n}$. Then
$\mathcal{G}(k)$ is the set of all potentials $f \in \mathcal{R}(k)$
such that, all periodic orbits with positive period in the
energy level $(H+f)^{-1}(k)$ are nondegenerate of all  orders
for $H+f$.
\begin{lemma}\label{NondegeneracyLemma}\textnormal{(Nondegeneracy Lemma)}
Given $k \in \mathbb{R}$ and $0< c \in \mathbb{R}$, the
set $\mathcal{G}_{k}^{c,c}$ is open and dense in
$C^{\infty}(M;\mathbb{R})$.
\end{lemma}
%###################################################################################
If $\mathcal{G}_{k}$ is generic, then  generically in $L$, the energy
level $k$ is regular and all periodic orbits in this 
level  are nondegenerate of all  orders for $H+f$. 
Thus, we must to prove that $\mathcal{G}_{k}^{c,c}$ is open in
 $C^{\infty}(M;\mathbb{R})$ , $\forall c \in \mathbb{R}_{+}$ and
dense in $\mathcal{R}(k)$, since Theorem~\ref{PotRegAbertoDenso},
implies that $\mathcal{R}(k)$ is dense in  $C^{\infty}(M;\mathbb{R})$. 
The proof of this lemma requires a sequence of technical constructions.
%###################################################################################
\begin{lemma}\label{MinimoPerioLema}
Given $k \in \mathbb{R}$ and $f_{0} \in \mathcal{R}(k)$ there exists a neighborhood, $\mathcal{U}$,
of $f_{0}$ in $C^{\infty}(M;\mathbb{R})$ and $0< \alpha:=\alpha(\mathcal{U},f_{0})$
such that, for all $f \in \mathcal{U}$,  the period of all periodic orbits of $H+f$, in the level $(H+f)^{-1}(k)$ , is bounded below by $\alpha$.
\end{lemma}
\begin{proof}
If we suppose that our claiming  is false, we get the existence of sequences,
$\mathcal{U} \ni f_{n} \to f_{0}$,  $ T_{n} >
0$ with $ T_{n} \to  0$ and $\vartheta_{n} \in (H+f_{n})^{-1}(k)$
such that $\psi_{T_{n}}^{H+f_{n}}(\vartheta_{n})=\vartheta_{n}$.

From Lemma~\ref{compacidade}, we can choose a subsequence such that
$$d_{T^{*}M}(\psi_{t}^{H+f_{n}}(\vartheta_{0}),
\vartheta_{0})=0, \, \forall t>0$$, that is, $\vartheta_{0} \in
(H+f_{0})^{-1}(k)$ which is a fix point, contradicting the fact of $f_{0} \in \mathcal{R}(k)$. 
\end{proof}
%###################################################################################
\begin{lemma}\label{openlema}
Given $k \in \mathbb{R}$, $a,b \in \mathbb{R}$ with $0 < a \leq b < \infty$, the set
$\mathcal{G}_{k}^{a,b}$ is open in  $C^{\infty}(M;\mathbb{R})$.
\end{lemma}
\begin{proof}
If, $\mathcal{G}_{k}^{a,b} \neq
\varnothing$, take $f_{0} \in
\mathcal{G}_{k}^{a,b}$. If $f_{0}$ is not an interior point we get the existence of a
sequence  $f_{n} \to f_{0}$ where
$f_{n} \not\in \mathcal{G}_{k}^{a,b}$. Therefore, there exists $\vartheta_{n} \in
(H+f_{n})^{-1}(k)$, $ T_{n}=T_{min}(\vartheta_{n}) \in
(0;a]$ and natural numbers $\ell_{n} \geq 1$ such that,
$\ell_{n}T_{n} \leq b$, $\psi_{\ell_{n}T_{n}}^{H+f_{n}}(\vartheta_{n})=
\vartheta_{n}$  and  $d_{\vartheta_{n}}
\psi_{\ell_{n}T_{n}}^{H+f_{n}}$ do not have 1 as eigenvalue
with algebraic multiplicity bigger than 1.
Consider $\mathcal{U}_{0}$ and $0 <
\alpha:=\alpha(\mathcal{U}_{0},f_{0}) <a$ as in
Lemma~\ref{MinimoPerioLema}. Choosing a subsequence we can assume that
$f_{n} \in \mathcal{U}_{0}$ and therefore $ T_{n} \in [\alpha;a]$,
with $\vartheta_{n} \to  \vartheta_{0} \in \varepsilon^{k}_{f_{0}}$,\, $T_{n} \to
T_{0}$,\, $\ell_{n} =  \ell_{0}$,\, $0< \alpha \leq T_{0} \leq a$\,
and \, $\ell_{0}T_{0} \leq b$.   Then $\psi_{\ell_{0}T_{0}}^{L+f_{0}}(\vartheta_{0})=\vartheta_{0}$, and
$d_{\vartheta_{0}}\psi_{\ell_{0}T_{0}}^{H+f_{0}}$ has 1 as eigenvalue
with algebraic multiplicity bigger than 1, that is, $\vartheta_{0}$
is a periodic orbit with minimal period $\leq$ $a$, degenerate of order
$\ell_{0} \leq \frac{b}{T_{0}}$,  contradicting the fact of
$f_{0} \in  \mathcal{G}_{k}^{a,b}$.
\end{proof}

%##############################################################################

In order to  prove that $\mathcal{G}_{k}^{c,c}$ is dense in
$C^{\infty}(M;\mathbb{R})$ , $\forall c \in \mathbb{R}_{+}$, we observe that is enough  show that, $\mathcal{G}_{k}^{c,c}$ is dense in
$\mathcal{R}(k)$. So, we can reduce this proof to a local approach.
More precisely, the claim is a direct consequence of the following
Reduction Lemma, whose proof we will  present in the Section~\ref{ReducLema}.
\begin{lemma}\label{ReduLocalDaPrimParte} \textnormal{(Reduction Lemma)} 
For each  $c \in \mathbb{R}_{+}$, and any $f_{0} \in \mathcal{R}(k)$,
there exists an open neighborhood $\mathcal{U}_{f_{0}}$ of $f_0$, such that,
$\mathcal{G}_{k}^{c,c} \cap \mathcal{U}_{f_{0}}$ is dense in $\mathcal{U}_{f_{0}}$.
\end{lemma}
Thus the Nondegeneracy Lemma is proven.\\

%###################################################################################
{\large \textbf{Proof of  the Kupka-Smale Theorem }}\\
%###################################################################################

   Consider a periodic orbit $\gamma=
\{\phi_{t}^{L}(\theta_{0}), \, 0 \leq t \leq T \} \subseteq
H^{-1}(k)$, in $T^{*}M$, where $H$ is the Hamiltonian 
associated to $L$ by the Legendre transform.
We will say that this orbit is hyperbolic if the Poincar\'e map
associated does not have eigenvalue of norm 1. It is clear that the
hyperbolicity implies in the nondegeneracy property. The
converse is not true. There exists  nondegenerate
orbits such that all eigenvalues has norm 1. Such orbits will be
called elliptic orbits.
We define the  strong stable and strong unstable manifolds, of
$\gamma$ in $\theta_{0}=\gamma(0)$, as
$$\displaystyle W^{ss}(\theta_{0})=\{ \theta \in H^{-1}(k) \mid
\lim_{t \to +\infty} d(\phi_{t}^{L}(\theta_{0}),\phi_{t}^{L}(\theta))=0 \}$$
and
$$\displaystyle W^{us}(\theta_{0})=\{ \theta \in H^{-1}(k) \mid
\lim_{t \to -\infty} d(\phi_{t}^{L}(\theta_{0}),\phi_{t}^{L}(\theta))=0 \}.$$

Respectively we define the stable and unstable manifolds (weak) of
$\gamma$ as $$\displaystyle W^{s}(\gamma)=\bigcup_{t \in
\mathbb{R}}\phi_{t}^{L}(W^{ss}(\theta_{0}))  \text{ and }  W^{u}(\gamma)=\bigcup_{t \in
\mathbb{R}}\phi_{t}^{L}(W^{us}(\theta_{0})).$$
From the general theory of the Lagrangians systems we know that,
$W^{s}(\gamma),$ $ W^{u}(\gamma) \subset H^{-1}(k)$ are Lagrangians submanifolds
of $TM$, with the symplectic twist form, given by
$\omega(\xi, \zeta)= \langle (\xi_{h}, \xi_{v})^{*}, J (\zeta_{h},
\zeta_{v}) \rangle $
in local coordinates.

A point $\theta \in  H^{-1}(k)$,  is heteroclinic if  $\theta
\in W^{s}(\gamma_{1}) \cap W^{u}(\gamma_{2})$. where $\gamma_{1}, \gamma_{2} \subset
H^{-1}(k)$ are hyperbolic periodic orbits. Additionally,
if $T_{\theta}W^{s}(\gamma_{1})$ $ + T_{\theta}W^{s}(\gamma_{2}) =
T_{\theta}H^{-1}(k)$, that is, if  $W^{s}(\gamma_{1}) \pitchfork_{\theta}
W^{u}(\gamma_{2})$  then $\theta$  will be called a transversal heteroclinic point, same
thing for homoclinics.

A fundamental domain for $W^{s}(\gamma)$ (or $W^{u}(\gamma)$)
is a compact subset $\mathcal{D} \subset W^{s}(\gamma)$,
such that, all orbits in $W^{s}(\gamma)$  intercepts
$\mathcal{D}$ in one point, at least. One can show that there exists fundamental
domains arbitrarily small and arbitrarily close to $\gamma$.
Fixed $a >0$ we define the local stable and local unstable
submanifolds of $\gamma$ as being
$$W^{s}_{a}(\gamma)=\{ \theta \in W^{ss}(\gamma) \mid
d_{W^{ss}(\gamma)}(\theta, \gamma) < a \}$$
and
$$W^{u}_{a}(\gamma)=\{ \theta \in W^{us}(\gamma) \mid
d_{W^{us}(\gamma)}(\theta, \gamma) < a \}.$$
They are Lagrangians submanifolds of $TM$.

In order to prove the  Kupka-Smale Theorem , we define
$\mathcal{K}^{a}_{k}=\{ f \in \mathcal{G}^{a,a}_{k} \mid
\forall \gamma_{1}, \gamma_{2} \subset  (H+f)^{-1}(k),$ hyperbolic periodic orbits\\
for $L+f$, with period $\leq a$ we have $W^{s}_{a}(\gamma_{1}) \pitchfork
W^{u}_{a}(\gamma_{2}) \}$   and
$\displaystyle \mathcal{K}(k)=\bigcap_{n \in \mathbb{N}} \mathcal{K}^{n}_{k}$.
It is clear that the properties (i), (ii) and (iii) of the Kupka-Smale Theorem   
are valid for all $f \in \mathcal{K}(k)$.
Thus, in order to prove the  Kupka-Smale Theorem , we must to show that $\mathcal{K}(k)$
is generic, or equivalent, that each, $\mathcal{K}^{n}_{k}$
is an open and dense set (in $C^{\infty}$ topology).
Since, the local stable and unstable  manifolds depends $C^{1}$ continuously 
on compact parts, of the Lagrangian field,
we get the openness of $\mathcal{K}^{n}_{k} $,
because the transversality is an open property.

The next lemma can be found  in Paternain~\cite{PaternainGeodFlows}, Proposition
2.11, Pg.34, for the geodesic case, but here we present a Lagrangian
version.

\begin{lemma} \label{PropTwistDoFibradoVertical} \textnormal{(Twist Property of the vertical
bundle)} Let $L$  be a smooth, convex and superlinear, Lagrangian
in $M$, $\theta \in TM$ and $F \subset T_{\theta} TM$ an Lagrangian
subspace for the twist form in $T^{*}M$. Then, the set,
$$\mathcal{Z}_{F}=\{ t \in \mathbb{R} \mid d_{\theta}\phi_{t}^{L}(E)
\cap V(\phi_{t}^{L}(\theta))
\neq \varnothing \}$$ is discrete, where $V$ is the vertical bundle in $M$.
\end{lemma}

The next lemma allow us to make a local perturbation of a potential
$f$ in such way that the correspondent stable and unstable manifolds
become transversal in a certain heteroclinic point $\theta$. The
density of $\mathcal{K}^{a}_{k}$ follows from
Lemma~\ref{DensidadeTransvVariedEstInst}.
%###################################################################################
\begin{lemma}\label{PertTransvVariedEstInst}
Let $L$ be a Lagrangian, and $f \in C^{\infty}(M,\mathbb{R})$. 
Given $\gamma_{1}, \gamma_{2} \subset
(H-f)^{-1}(k)$ hyperbolic periodic orbits with period $\leq a$ and
$\theta \in W_{a}^{u}(\gamma_{2})$, such that, the canonic
projection $\pi \mid_{W_{a}^{u}(\gamma_{2})}$ is a local
diffeomorphism in $\theta$ and $U, V$, are neighborhoods of 
$\theta$ in $TM$ such that $\theta \in V \subset
\bar{V} \subset U$.
Then, there exists $\bar{f} \in C^{\infty}(M,\mathbb{R})$, such that,\\
i) $\bar{f} $ is $ C^{\infty}$ close to $f$;\\
ii) $\supp(  f - \bar{f}) \subset \pi (U)$; \\
iii) $\gamma_{1}, \gamma_{2} \subset  (H-\bar{f})^{-1}(k)$ are
hyperbolic periodic orbits to $\bar{f}$, with the same period as to $f$; \\
iv) The connected component of $W_{a}^{u}(\gamma_{2}) \cap V$ that
contains $\theta$ is transversal to $W^{s}(\gamma_{1})$.
\end{lemma}
\begin{proof}
Initially we consider the Hamiltonian $H-f$ associated to the Lagrangian
$L+f$ by the Legendre transform $\mathcal{L}$:
$$ H-f(x,p)=\sup_{v \in T_{x}M} \{ p(v)- (L+f)(x,v)\} $$
with the canonic symplectic form of $T^{*}M$, $\omega = \sum dx_{i} \wedge dp_{i}.$

We know, from the general theory of the Hamiltonian systems, that
$\gamma_{1}, \gamma_{2}$ are in correspondence, by Legendre
transform with hyperbolic periodic orbits of same period,
$\tilde{\gamma_{1}}, \tilde{\gamma_{2}} \subset (H-f)^{-1}(k)$ for
the Hamiltonian flow $\psi_{t}^{H-f}$. Consider
$\tilde{W}^{u}(\tilde{\gamma_{2}})$ and
$\tilde{W}^{s}(\tilde{\gamma_{1}})$, respectively, the invariant
submanifolds, they will be Lagrangian submanifolds of $T^{*}M$, and
$\vartheta= \mathcal{L}(\theta) \in
\tilde{W}^{u}(\tilde{\gamma_{2}})$.
If $\pi:TM \to M$ and $\pi^{*}:T^{*}M \to M$ are the
canonic projections, then
$d_{\vartheta}\pi^{*}= d_{\theta}\pi  \circ (d_{\theta}\mathcal{L})^{-1}$.
Therefore the canonic projection $\pi^{*}
\mid_{\tilde{W}^{u}(\tilde{\gamma_{2}})}$ is a local diffeomorphism
in $\vartheta$. Moreover, $X^{H-f}(\vartheta)=d_{\theta}\mathcal{L}
\circ X^{L+f}(\theta) \neq 0$. Thus we can prove the lemma in the
Hamiltonian setting.
By \cite{GonzaloPaternGenGeodFlowPositEntropy}, Lemma A3,  we can
find a neighborhood $ U$ of $\vartheta$, and $V \subset U$, such
that, $ V \subset \bar{V} \subset U$, and a Lagrangian submanifold,
$\mathcal{N}$, $C^{\infty}$ close to
$\tilde{W}^{u}(\tilde{\gamma_{2}})$, satisfying the following conditions\\
1) $\vartheta \in \{ U \backslash \bar{V} \}$\\
2) $\mathcal{N} \cap \{ U \backslash \bar{V} \} =
\tilde{W}^{u}(\tilde{\gamma_{2}}) \cap \{ U \backslash \bar{V} \}
\subset (H - f)^{-1}(k)$;\\
3) $\mathcal{N} \cap \bar{V}  \pitchfork
\tilde{W}^{s}(\tilde{\gamma_{1}})  \cap \bar{V}$.

As, $\mathcal{N}$ is  $C^{\infty}$ close to
$\tilde{W}^{u}(\tilde{\gamma_{2}})$, we have that the canonic
projection $\pi^{*} \mid_{\mathcal{N}}$ is a local
diffeomorphism in $\vartheta$. If $U$ is
small enough, then
$\mathcal{N} \cap U =\{ (x, p(x)) \mid x \in \pi^{*}(u)\}$ 
that is, $\mathcal{N}\mid_{U}$ is a $C^{\infty}(M,\mathbb{R})$
graph.  We define the following potential, $\bar{f} \in
C^{\infty}(M,\mathbb{R})$,
$$
\bar{f}(x)= \left \{
\begin{array}{clcr}
f(x) \quad \, \quad \quad \quad \quad if \; x \in \pi^{*}(U)^{C}\\
H(x,p(x))-k \quad \, if \; x \in \pi^{*}(U) \quad _{}
\end{array}
\right.
$$

Observe that, $\supp(  f - \bar{f}) \subset \pi^{*}(U)$ and
$\vartheta \not\in \supp( f - \bar{f})$, moreover, choosing $U$ small
enough, we will have that $\pi^{*}(U) \cap \{ \tilde{\gamma_{1}},
\tilde{\gamma_{2}} \} =\varnothing$ and therefore
$\tilde{\gamma_{1}}, \tilde{\gamma_{2}}$ still, hyperbolic periodic
orbits of same period for the Hamiltonian flow
$\psi_{t}^{H-\bar{f}}$, contained in $(H-\bar{f})^{-1}(k)$.
We denote $\bar{W}^{u}(\tilde{\gamma_{2}})$  and
$\bar{W}^{s}(\tilde{\gamma_{1}})$, the invariant  manifolds for the
new flow $\psi_{t}^{H-\bar{f}}$.
Clearly $(H-\bar{f})(\mathcal{N})=k$. By
\cite{GonzaloPaternGenGeodFlowPositEntropy}, Lemma A1, we have that
$\mathcal{N}$ is $\psi_{t}^{H-\bar{f}}$ invariant. Since,
$\bar{W}^{u}(\tilde{\gamma_{2}})$ depends only of the negative times
and, the connected component of $\bar{W}^{u}(\tilde{\gamma_{2}})
\cap U$ that contains $\vartheta$ and $\mathcal{N}$ are coincident
in a neighborhood of $ \tilde{\gamma_{2}}$ disjoint of  $\supp( f -
\bar{f})$, we have 
$\mathcal{N}=\bar{W}^{u}(\tilde{\gamma_{2}})$.
On the other hand, as $\bar{W}^{s}(\tilde{\gamma_{1}})$ depends only
of the positive times and $f=\bar{f}$ in $\{ U \backslash \bar{V}
\}$, we have $\bar{W}^{s}(\tilde{\gamma_{1}}) =
\tilde{W}^{s}(\tilde{\gamma_{1}})$. Since $\mathcal{N} \cap \bar{V}
\pitchfork \tilde{W}^{s}(\tilde{\gamma_{1}})  \cap \bar{V}$,  we have
$\bar{W}^{u}(\tilde{\gamma_{2}}) \cap \bar{V} \pitchfork
\tilde{W}^{s}(\tilde{\gamma_{1}})  \cap \bar{V}.$ From the
initial considerations we choose $L+ \bar{f} $. The lemma is proven.
\end{proof}

%###################################################################################
\begin{lemma}\label{DensidadeTransvVariedEstInst}
The set $\mathcal{K}^{n}_{k}$, is dense in
$C^{\infty}(M,\mathbb{R})$, for all $n \in \mathbb{N}$.
\end{lemma}
\begin{proof}
Take $f_{0} \in C^{\infty}(M,\mathbb{R})$, by the Nondegeneracy Lemma
we can  find $f_{0}$ arbitrarily  close to $f' \in
\mathcal{G}_{k}^{n,n}$, which  is open and dense. Thus, is enough 
to find $f$ arbitrarily close to $f'$, such that, for any
$\gamma_{1}, \gamma_{2} \subset (H-f)^{-1}(k)$, hyperbolic periodic
orbits of period $\leq n$, is valid  $W_{n}^{s}( \gamma_{1})
\pitchfork W_{n}^{u}( \gamma_{2})$. Then, $f \in
\mathcal{K}^{n}_{k}$ and $f$ is arbitrarily close to $f_{0}$.
Given $\gamma_{1}, \gamma_{2} \subset (H+f')^{-1}(k)$
hyperbolic periodic orbits of period $\leq n$, in order to conclude
that $W_{n}^{s}( \gamma_{1}) \pitchfork W_{n}^{u}( \gamma_{2})$ we 
should to prove that $W_{n}^{s}( \gamma_{1})
\pitchfork_{\mathcal{D}} W_{n}^{u}( \gamma_{2})$ where $\mathcal{D}$
is a fundamental domain of $W^{u}( \gamma_{2}) $ , because if
$W^{s}( \gamma_{1}) \pitchfork_{\theta} W^{u}( \gamma_{2})$ then
$W^{s}( \gamma_{1}) \pitchfork_{\phi_{t}^{L+f'}(\theta)} W^{u}(
\gamma_{2}), \, \forall t$.
Take $\mathcal{D}$ a fundamental domain of $W^{u}( \gamma_{2}) $ and
$\theta \in \mathcal{D}$. By the inverse function theorem we know
that $\pi|_{W^{u}( \gamma_{2})}$ is a local diffeomorphism in
$\theta$ if, and only if, $T_{\theta}W^{u}( \gamma_{2}) \cap Ker
d_{\theta} \pi ={0}$. As $W^{u}( \gamma_{2})$ is a Lagrangian
submanifold we have, from Lemma~\ref{PropTwistDoFibradoVertical},
that
$\{ t \in \mathbb{R} \mid d_{\theta}\phi_{t}^{L}(T_{\theta}W^{u}( \gamma_{2}))
\cap Ker \, d_{\phi_{t}^{L+f'}(\theta)} \pi \neq \varnothing \},$ is
discrete. Then there exists $t(\theta)$ arbitrarily close to $0$, such
that, $\pi|_{W^{u}( \gamma_{2})}$ is a local diffeomorphism in
$\tilde{\theta}=\phi_{t(\theta)}^{L+f'}(\theta)$.
As, $f' \in \mathcal{G}_{k}^{n,n}$, we can choose, $t(\theta)$, such
that, $\pi(\tilde{\theta})$ does not intercept any periodic orbit of
period  $\leq n$. Fix a  neighborhood $U$, of $\tilde{\theta}$,  arbitrarily small,
such that, $\pi(U)$ does not intercept any
periodic orbit of period  $\leq n$. Taking $V$, a neighborhood of
$\tilde{\theta} $, such that,  $V \subset \bar{V} \subset U$, from
Lemma~\ref{PertTransvVariedEstInst}, we can find $f_{1}=f'$ in
$\pi(U)^{C}$, such that, the connected component of
$W_{n}^{u}( \gamma_{2}) \cap \bar{V}$ (to the new flow)
contain  $\tilde{\theta}$, and is transversal to $W^{s}( \gamma_{1})$ 
(to the new flow). Taking $\bar{V}_{1} =
\phi_{t'}^{L+f_{1}}(\bar{V})$ we will have that $W_{n}^{u}(
\gamma_{2}) \pitchfork W^{s}( \gamma_{1})$.
We can cover the fundamental domain $\mathcal{D}$ with a finite number
of neighborhoods like $\bar{V}_{1}$, that is, $W_{1},..., W_{s}$.
Since the transversality is an open condition and the local
stable (unstable) manifold depends continuously on compact parts, we
can choose successively $W_{i+1}$ such that the transversality in
$W_{j}, \; j \leq i$, is preserved. Thus,  
$W_{n}^{u}( \gamma_{2}) \pitchfork W_{n}^{s}( \gamma_{1})$.
\end{proof}

%###################################################################################
\section{Proof of the Reduction Lemma} \label{ReducLema} 
%###################################################################################

For the proof of the Reduction Lemma
(Lemma~\ref{ReduLocalDaPrimParte}) we will use an induction
method similar to the one used by Anosov~\cite{Anos}, using transversality arguments
as described in Abraham \cite{Ab} and \cite{AbRb}. In this way, we remember a usefull theorem, the Parametric Transversality Theorem of Abraham.

Remember that, if $\mathcal{X}$ is a topological space. A subset $\mathcal{R} \subseteq
\mathcal{X}$ is said generic if $\mathcal{R}$ is a countable
intersection of open and dense sets. The space
$\mathcal{X}$ will be a Baire Space if all generic subsets are dense.  For additional results and definitions of Differential Topology, see~\cite{AbRb}, \cite{AbMardRat} or \cite{KbMg}.
\begin{theorem}\label{abraham}(\cite{AbRb}, pg. 48,  Abraham's Parametric Transversality Theorem)
Consider $\mathbb{X}$ a submanifold finite dimensional (with
boun\-da\-ry or boundaryless), $\mathbb{Y}$ a boun\-da\-ry\-less manifold and $S
\subseteq \mathbb{Y}$ a submanifold with finite codimension.
Consider $\mathcal{B}$ boundaryless manifold, $\rho : \mathcal{B}
\rightarrow C^{\infty}(\mathbb{X};\mathbb{Y})$ a smooth
representation and your evaluation $ev_{\rho} : \mathcal{B} \times
\mathbb{X} \rightarrow \mathbb{Y}$. If $\mathbb{X}$ and
$\mathcal{B}$ are Baire spaces and $ev_{\rho} \pitchfork S $ then
the set $\mathcal{R}=\{ \varphi \in \mathcal{B} \mid \rho_{\varphi}
\pitchfork S \}$ is a generic subset (and obviously dense) of
$\mathcal{B}$.
\end{theorem}

%###################################################################################

Given a Hamiltonian $H$ we define the normal field associated,
as being the gradient field, $Y^{H}=\nabla H= H_{x} \, \frac{\partial \,}{\partial x} + H_{p} \, \frac{\partial \,}{\partial p}$ in $T^{*}M$. Observe that $JY^{H}=X^{H}$, where  $J$ is the canonic sympletic matrix. We denote $\psi_{s}^{H^{\perp}}: T^{*}M \times (-\varepsilon,\varepsilon) \rightarrow T^{*}M$ the flow in $T^{*}M$ generated by the normal field.
  Let us briefly describe  the properties of the normal field. Initially
  observe that $\omega_{\vartheta} (Y^{H}, X^{H})= H_x^{2} +
  H_p^{2}, \; \forall \vartheta \in T^{*}M$. If $X^{H}(\vartheta) \neq
  0$ then $0 \neq Y^{H}(\vartheta) \not\in T_{\vartheta}H^{-1}(k)$ where
  $k=H(\vartheta)$, that is $Y^{H}$  points to the outside of the energy level.
  From the compactness of the energy level $H^{-1}(k)$
  we have that the flow of the normal field, restricted
   to $H^{-1}(k)$ is defined in  $H^{-1}(k) \times
  (-\varepsilon(H),\varepsilon(H))$ where  $\varepsilon(H) >0$ is uniformly defined in
  $H^{-1}(k)$. Then the flow of the normal
  field is defined in a neighborhood of the energy level $H^{-1}(k)$.
  The action of the differential of the normal flow through an orbit is given by
$$
\left \{
\begin{array}{clcr}
\dot{Z}^{H}(s) = \mathcal{H}(\gamma(s)) Z^{H}(s) \\
Z^{H}(0)=Y^{H}(\vartheta) \quad  \quad \, \, \; \; \; ^{}
\end{array}
\right .
$$
where $\mathcal{H}$ is the hessian matrix of $H$.  The main property  of $Y^{H}$  is to establish a sympletic decomposition of $T_{\vartheta}T^{*}M$ given by the next proposition.

\begin{proposition}\label{BaseCampoNormal} Consider the normal field
$Y^{H}$ associated to $H$ in the regular energy level, $H^{-1}(k)$.
Then, for each $\vartheta \in H^{-1}(k) $, periodic of period $T
>0$, there exists a symplectic base $\{
u_{1},...,u_{n},u_{1}^{*},...,u_{n}^{*} \}$ of
$T_{\vartheta}T^{*}M $ verifying\\
(i) $u_{1}= X^{H} $ and  $u_{1}^{*}=-\frac{1}{H_x^{2} +
  H_p^{2}} Y^{H} $;\\
(ii) $\mathcal{W}_{1}= \langle u_{1}, u_{1}^{*} \rangle ^{\perp}
\subset T_{\vartheta} H^{-1}(k) $ in particular, $T_{\vartheta}
H^{-1}(k) = \langle u_1 \rangle \oplus \mathcal{W}_{1}$;\\
(iii) If $\Sigma \subset H^{-1}(k)$ is a section transversal to the
flow, such that, $T_{\vartheta} \Sigma= \mathcal{W}_{1}$, 
we have that  $d_{\vartheta} P(\Sigma,\vartheta)
\mathcal{W}_{1} \subseteq \mathcal{W}_{1} $;\\
(iv) If $T=m T_{min}(\vartheta)$, then
$d_{\vartheta} \psi_{T}^{H} u_{1}=u_{1}$ and $d_{\vartheta}
\psi_{T}^{H} u_{1}^{*}=c u_{1} + u_{1}^{*} +\xi$, $\xi \in
\mathcal{W}_{1}$.
In particular, $(d_{\vartheta} \psi_{T}^{H} -Id)
(T_{\vartheta}T^{*}M) \subseteq \langle u_{1} \rangle \oplus
\mathcal{W}_{1}=T_{\vartheta} H^{-1}(k) $.\\
(v) There exists, $\varepsilon >0$ uniform in $\vartheta \in H^{-1}(k)$, such
that, the map $e_{\vartheta} : (-\varepsilon, \varepsilon) \to~\mathbb{R}$ given by
$e_{\vartheta}(s)=H \circ \psi_{s}^{H^{\perp}}(\vartheta)$ is
injective with $e_{\vartheta}(0)=k$.
\end{proposition}

Using the normal field we are able to construct an representation, in order  to apply the Parametric Transversality Theorem.

\begin{proposition}\label{repres}
 Given $k \in \mathbb{R}$, $0<a<b< + \infty$, and $f_{0} \in \mathcal{R}(k)$,
 consider the normal field $Y^{H+f_0}$ as described before,
 $\varepsilon=\varepsilon(H+f_0) >0$ as in the Proposition~\ref{BaseCampoNormal},~(v),
 and the sets, $\mathcal{U}_{f_0} \subset \mathcal{R}(k)$ a $C^{\infty}$ neighborhood
  of $f_0$, $\alpha=\alpha(\mathcal{U}_{f_0}) >0$
 as in Lemma~\ref{MinimoPerioLema}, $\mathbb{X}= T^{*}M \times
 (a,b) \times (-\varepsilon, \varepsilon)$ and $\mathbb{Y}= T^{*}M \times T^{*}M \times
 \mathbb{R}$. 
 Then the map  $\rho : \mathcal{U}_{f_0} \rightarrow
 C^{\infty}(\mathbb{X};\mathbb{Y})$ given by
 $\rho(f):=\rho_{f}$, where
 \begin{center}$\rho_{f}(\vartheta, t, s)=(\psi_{s}^{H+f^{\perp}}
 (\vartheta), \psi_{t}^{H+f}(\vartheta),
 (H+f)(\vartheta)-k)$ \end{center}
 is an injective representation (see~\cite{Ab}~or~\cite{AbRb}).
\end{proposition}

\begin{proof}
Initially  we point out that $\rho$
is well defined, therefore $\mathbb{Y}$  has the structure of a product manifold. Writing
$\rho_{f}=(\rho_{f}^{1},\rho_{f}^{2},\rho_{f}^{3})$, where
 $\rho_{f}^{1}(\vartheta, t,s)=\psi_{s}^{H+f^{\perp}}(\vartheta) $, 
$\rho_{f}^{2}(\vartheta, t,s) = \psi_{t}^{L+f}(\vartheta) $,
$\rho_{f}^{3}(\vartheta, t,s) = (H+f)(\vartheta)-k$,
we can see that each coordinate is a smooth  function. Thus
$\rho_{f}\in C^{\infty}(\mathbb{X};\mathbb{Y})$.
Observe that $\rho$ is injective. Indeed, if $\rho_{f_{1}}=\rho_{f_{2}}$
then  $(H+f_{1})(\vartheta)-k =
(H+f_{2})(\vartheta)-k$, for all $\vartheta \in T^{*}M$, so $f_{1}(x) = f_{2}(x)
$, for all $x \in M$.  Thus $f_{1} = f_{2}$.
We must to verify that $ev_{\rho} : \mathcal{U}_{f_0} \times \mathbb{X}
\to \mathbb{Y}$  is smooth. Since $\mathcal{U}_{f_0} \times
\mathbb{X}$ have the  structure  of product
manifold , we can write
$$d_{(f,x)} ev_{\rho}:= \frac{\partial ev_{\rho}}{\partial f} (f,x)
+ \frac{\partial ev_{\rho}}{\partial x} (f,x), \; \forall
(f,x)\in \mathcal{U}_{f_0} \times \mathbb{X}$$
It is clear that $\frac{\partial ev_{\rho}}{\partial x}
(f,x)$ is always defined as$ \frac{\partial ev_{\rho}}{\partial x} (f,x)=d_{x}\rho_{f}$.
More precisely, given $(\xi,\dot{t},\dot{s}) \in T_{x=(\vartheta,T,S)}\mathbb{X}$
we have
$$\frac{\partial ev_{\rho}}{\partial x}
(f,x)(\xi,\dot{t},\dot{s})=d_{x}\rho_{f}(\xi,\dot{t},\dot{s}) =\frac{d \,}{dr}\rho_{f}(\vartheta(r), t(r), s(r)) \mid_{r=0}=$$
$$=\frac{d \,}{dr} (\psi_{s(r)}^{H+f^{\perp}}(\vartheta(r)), \psi_{t(r)}^{H+f}(\vartheta(r)),
(H+f)(\vartheta(r))-k)\mid_{r=0}=$$
$$( d_{\vartheta} \psi_{S}^{H+f^{\perp}}(\xi) +
\dot{s} Y^{H+f}(\psi_{S}^{H+f^{\perp}}(\vartheta)), d_{\vartheta}\psi_{T}^{H+f}(\xi) +
\dot{t} X^{H+f}(\psi_{S}^{H+f}(\vartheta)),$$ $$ d_{\vartheta}(H+f)(\xi) ).$$
Observe that, if $S=0$ and $\psi_{S}^{H+f}(\vartheta)=\vartheta$,
then  $$\frac{\partial ev_{\rho}}{\partial x}
(f,x)(\xi,\dot{t},\dot{s})=( \xi + \dot{s} Y^{H+f}(\vartheta), d_{\vartheta}\psi_{T}^{H+f}(\xi) +
\dot{t} X^{H+f}(\vartheta),$$ $$ d_{\vartheta}(H+f)(\xi) ).$$ 
However we must to show that $\frac{\partial ev_{\rho}}{\partial f}
(f,x)$ is always defined. By the structure of  $C^{\infty}(M;\mathbb{R})$
we know that this fact  is equivalent to show that there exists
\begin{center} $\frac{d \,}{dr} \psi_{S}^{H + f + r h ^{\perp}}(\vartheta) \mid_{r=0}$ ,
$ \frac{d \,}{dr} \psi_{T}^{H + f + r h }(\vartheta)
\mid_{r=0}$ and $\frac{d \,}{dr} (H + f + r h ) (\vartheta) \mid_{r=0}$
\end{center}
for any $h \in C^{\infty}(M;\mathbb{R})$  and
$x=(\vartheta,T, S) \in \mathbb{X}$.  From some straightforward  
calculations (see~\cite{PhilHart}, pg.46), we get\\
$\frac{d \,}{dr} (H + f + r h ) (\vartheta) \mid_{r=0}=h\circ\pi(\vartheta)$,\\
$ \frac{d \,}{dr} \psi_{T}^{H + f + r h }(\vartheta)
\mid_{r=0}=Z_{h}(T)= d_{\vartheta}\psi_{T}^{H + f}
\int_{0}^{T} (d_{\vartheta}\psi_{t}^{H + f })^{-1}
b_{h}(t)dt$ and\\
$ \frac{d \,}{dr} \psi_{S}^{H + f + r h ^{\perp} }(\vartheta)
\mid_{r=0}=Z^{h}(S)= d_{\vartheta}\psi_{S}^{H + f^{\perp} }
\int_{0}^{S} (d_{\vartheta}\psi_{s}^{H + f^{\perp}  })^{-1}
b^{h}(s)ds$.\\
Thus  $\displaystyle  \frac{\partial ev_{\rho}}{\partial f} (f,x)
(h)=$
 $$ ( d_{\vartheta}\psi_{S}^{H + f^{\perp} }
\int_{0}^{S} (d_{\vartheta}\psi_{s}^{H + f^{\perp}  })^{-1}
b^{h}(s)ds \; ,   \; d_{\vartheta}\psi_{T}^{H + f}
\int_{0}^{T} (d_{\vartheta}\psi_{t}^{H + f })^{-1}
b_{h}(t)dt  \; ,$$ $$  \; h\circ\pi(\vartheta)).$$
   If $S=0$ and $\psi_{T}^{H+f}(\vartheta)=\vartheta$,
   then $$\displaystyle  \frac{\partial ev_{\rho}}{\partial f} (f,x)
(h)=( 0  ,  \; d_{\vartheta}\psi_{T}^{H + f}
\int_{0}^{T} (d_{\vartheta}\psi_{t}^{H + f })^{-1}
b_{h}(t)dt  \; ,  \; h \circ \pi(\vartheta)).$$
Thus $ev_{\rho}$ is smooth
and therefore $\rho$ is a representation.
\end{proof}

Define the null diagonal $\Delta_{0} \subseteq
\mathbb{Y}$ given by $\Delta_{0}=\{(\vartheta,\vartheta,0) \mid \vartheta \in T^{*}M\}$.
Combining, Propositions~\ref{BaseCampoNormal} and \ref{repres} we get

\begin{lemma}\label{PeriodPontoSaoNaoDegenTrans}
With the same notations of the Proposition~\ref{repres}, we have that,
$\forall f \in \mathcal{U}_{f_0}$, with $T \in (a,b)$ and $S \in (-\varepsilon, \varepsilon)$,
\begin{itemize}
\item[i)] If $\vartheta$ is a periodic orbit of positive period $T$ for $H+f$
in the level $(H+f)^{-1}(k)$ then,
$\rho_{f}(\vartheta, T, 0) \in \Delta_{0}$.
Reciprocally, if $\rho_{f}(\vartheta, T, S) \in \Delta_{0}$
then, $S=0$ and $\vartheta$ is a periodic orbit of positive period $T$
for $H+f$ in the level $(H+f)^{-1}(k)$.

\item[ii)] If $\vartheta$ is a periodic orbit of positive period, for $H+f$
in the level $(H+f)^{-1}(k)$. Then, $\vartheta$ is
nondegenerate of order $ m = \frac{T}{T_{min}(\vartheta)}$ if, and only if,
$\rho_{f} \pitchfork _{(\vartheta, T, 0)} \Delta_{0}$.
\end{itemize}
\end{lemma}

The next corollary it is an easy consequence of the Lemma~\ref{PeriodPontoSaoNaoDegenTrans}.
%##############################################################################

\begin{corollary}\label{PeriodNivelSaoNaoDegenTrans}
With the same notations of the Lemma~\ref{PeriodPontoSaoNaoDegenTrans} we have that,
given $ f \in \mathcal{U}_{f_0}$, all periodic orbits $\vartheta $, with positive period,
$T_{min}(\vartheta) \in (a,b)$, in $(H+f)^{-1}(k)$,
are nondegenerate for $H+f$ of order $m$, $\forall m \leq \frac{b}{T_{min}}$
if, and only if, $\rho_{f} \pitchfork \Delta_{0}$.
\end{corollary}

%###################################################################################
The previous corollary shows that, the nondegeneracy of the periodic orbits
of positive period in an interval $(a,b)$, for a given energy level $(H+f)^{-1}(k)$,
is equivalent to the transversality of the map $\rho_{f}$ in relation to the diagonal $\Delta_{0}$. 
The key element for the proof of the Lemma~\ref{ReduLocalDaPrimParte}  is
the nest lemma.
%###################################################################################
\begin{lemma}\label{evaluationtransverlema}
Consider  the representation $\rho$ as in the Proposition~\ref{repres} and its evaluation
in $\mathcal{U}_{f_0}$, that is, $ev : \mathcal{U}_{f_0} \times \mathbb{X} \rightarrow
\mathbb{Y}$, given by $ev(f, \vartheta, t , s)= \rho_{f}(\vartheta, t ,s)$.
Suppose that $ev(f , \vartheta, T , S) \in \Delta_{0}$ then,
\begin{itemize}
\item[i)] If $\vartheta$ is  nondegenerate of order $m=\frac{T}{T_{min}}$ for $H+f$
then $ev \pitchfork _{(f, \vartheta, T ,S)}
\Delta_{0}$;
\item[ii)] If $T = T_{min}(\vartheta)$ then, $ev \pitchfork _{(f, \vartheta, T, S )}
\Delta_{0}$.
\end{itemize}
\end{lemma}

\begin{proof}

i) We know that $ev(f, \vartheta, T ,S)=\rho_{f}(\vartheta, T ,S)$
therefore $\rho_{f}(\vartheta, T, S )\in  \Delta_{0}$, and  $S=0$.
If $\vartheta$ is nondegenerate of order $m=\frac{T}{T_{min}}$
for $H+f$, then, from the Lemma~\ref{PeriodPontoSaoNaoDegenTrans}~(ii),
$\rho_{f} \pitchfork _{(\vartheta, T, 0 )} \Delta_{0}$,
in particular  $ev \pitchfork _{(f, \vartheta, T, 0)}
\Delta_{0}$.\\

ii) As $ev(f, \vartheta, T , S) \in \Delta_{0}$ we must to show that
\begin{center}
$d_{(f, \vartheta, T , 0)} ev T_{(f, \vartheta, T, 0 )} (\mathcal{U}_{f_0} \times
\mathbb{X})
+ T_{(\vartheta , \vartheta, 0)} \Delta_{0}= T_{(\vartheta , \vartheta, 0)}\mathbb{Y} $.
\end{center}
Take any $(u,v,w ) \in T_{(\vartheta,\vartheta
,0)}\mathbb{Y}$, $(\zeta, \zeta, 0) \in T_{(\vartheta,\vartheta
,0)} \Delta_{0} $ and \\ $(h, \xi, \dot{t}, \dot{s}) \in T_{(f, \vartheta, T, 0 )}
(\mathcal{U}_{f_0} \times \mathbb{X})$. From Proposition~\ref{repres} we have that\\
$d_{(f,\vartheta, T , 0)} ev_{\rho}(h, \xi, \dot{t},
\dot{s})=$\\
$=( 0  ,  \; d_{\vartheta}\psi_{T}^{H + f}
\int_{0}^{T} (d_{\vartheta}\psi_{t}^{H + f })^{-1}
b_{h}(t)dt  \; ,  \; h\circ\pi(\vartheta))
+ \\
( \xi + \dot{s} Y^{H+f}(\vartheta), d_{\vartheta}\psi_{T}^{H+f}(\xi) +
\dot{t} X^{H+f}(\vartheta), d_{\vartheta}(H+f)(\xi) )$\\
$=( \xi + \dot{s} Y^{H+f}(\vartheta) , \;  d_{\vartheta}\psi_{T}^{H+f}(\xi) +
\dot{t} X^{H+f}(\vartheta)+ d_{\vartheta}\psi_{T}^{H + f}
\int_{0}^{T} (d_{\vartheta}\psi_{t}^{H + f })^{-1}
b_{h}(t)dt  \; , \;\\ h\circ\pi(\vartheta) +  d_{\vartheta}(H+f)(\xi) )$\\
Therefore $ev_{\rho} \pitchfork _{(f, \vartheta, T, 0 )}
\Delta_{0}$, if and only if, the system
\begin{center}
\(
\left \{
\begin{array}{clcr}
u=\xi + \dot{s} Y^{H+f}(\vartheta) +\zeta  \quad \quad \quad \quad \quad \quad \quad \quad
\quad \quad \quad \quad \quad \quad \quad \quad \quad \quad  \quad (1)\\
v=d_{\vartheta}\psi_{T}^{H+f}(\xi) +
\dot{t} X^{H+f}(\vartheta)+ d_{\vartheta}\psi_{T}^{H + f}
\int_{0}^{T} (d_{\vartheta}\psi_{t}^{H + f })^{-1}
b_{h}(t)dt + \zeta \; \; \, (2)\\
w=h\circ\pi(\vartheta) +d_{\vartheta}(H+f)(\xi)  \quad \quad \quad \quad
\quad \quad \quad \quad \quad \quad \quad \quad \quad \quad \quad   \quad(3)
\end{array}
\right.
\)
\end{center}
has a solution.Using the coordinates of the Proposition~\ref{BaseCampoNormal}
and  taking $\zeta=u -\xi - \dot{s} Y^{H+f}(\vartheta) $ we have that the equation (2)
restricted to the set of the solutions of (3),
$$\mathcal{V}_{w}=\left\{h \in C^{\infty}(M,\mathbb{R}), \, \xi = a X^{H+f}(\vartheta)
 + b_{0} Y^{H+f}(\vartheta) + U \; \right\} $$
where
$b_{0}=\frac{w-h\circ\pi(\vartheta)}{d_{\vartheta}(H+f)(Y^{H+f}(\vartheta))}$,
will have the expression\\
$(\dot{t}+b_{0}c+\tau_{0}) X^{H+f}(\vartheta)+ (c^{*}-  \dot{s}) Y^{H+f}(\vartheta) +
(d_{\vartheta} P(\Sigma, \vartheta)-Id)(U) +  b_{0}U_{0} + \\
d_{\vartheta}\psi_{T}^{H + f}\int_{0}^{T} (d_{\vartheta}\psi_{t}^{H + f })^{-1}
b_{h}(t)dt=  \tilde{a} X^{L+f}(\vartheta) + \tilde{b}  Y^{H+f}(\vartheta) +
\bar{U}$,  where
$$v-u= \tilde{a} X^{L+f}(\vartheta) + \tilde{b}  Y^{H+f}(\vartheta) + \bar{U},$$\
$$(d_{\vartheta}\psi_{T}^{H+f} -id)(Y^{H+f}(\vartheta))=c
X^{H+f}(\vartheta) + c^{*} Y^{H+f}(\vartheta) + U_{0}$$ and 
$$(d_{\vartheta}\psi_{T}^{H+f} -id)(U)=\tau_{0}
X^{H+f}(\vartheta) + (d_{\vartheta} P(\Sigma,
\vartheta)-Id)(U).$$
That is, the system  always has a solution, if the expression,
\begin{center}
$(\dot{t}+b_{0}c+\tau_{0}) X^{H+f}(\vartheta) + (c^{*}-  \dot{s}) Y^{H+f}(\vartheta) +
(d_{\vartheta} P(\Sigma, \vartheta)-Id)(U) +  b_{0}U_{0} +
d_{\vartheta}\psi_{T}^{H + f}\int_{0}^{T} (d_{\vartheta}\psi_{t}^{H + f })^{-1}
b_{h}(t)dt$
\end{center}
is surjective in $T_{\vartheta}T^{*}M$.
So, we must to show that $$ d_{\vartheta}\psi_{T}^{H + f}
\int_{0}^{T} (d_{\vartheta}\psi_{t}^{H + f })^{-1}
b_{h}(t)dt$$ generates  a $2n-2$ dimensional space complementary
to  the space generated by $ X^{H+f}(\vartheta)$ and $ Y^{H+f}(\vartheta) $,  in
$T_{\vartheta}T^{*}M$, which is the claim of the next lemma.
\end{proof}

\begin{lemma}\label{sobrejetividade}
With the same notations as in Lemma~\ref{evaluationtransverlema}, the map
$ \mathcal{B}: C^{\infty}(M ; \mathbb{R})$
$\to~T_{\vartheta}T^{*}M$,
\begin{center}
$\displaystyle \mathcal{B}(h)=- d_{\vartheta}\psi_{T}^{H + f}
\int_{0}^{T} (d_{\vartheta}\psi_{t}^{H + f })^{-1}
b_{h}(t)dt$
\end{center}
generates a space complementary to
$\langle X^{H+f}(\vartheta), Y^{H+f}(\vartheta) \rangle$.
\end{lemma}
\begin{proof}
In order to prove this claim is enough to restrict the map $\displaystyle
\mathcal{B}$ to a subspace chosen in $C^{\infty}(M ; \mathbb{R})$.
Consider $t_{0} \in (0,T)$,  $\varepsilon >0$, and
denote $\mathcal{A}_{t_{0}}$, the  subspace of the smooth functions
$$\mathcal{A}_{t_{0}}=\{\alpha: \mathbb{R} \to \mathbb{R}^{n-1} \mid
\alpha(t)=(a_{1}(t),...,a_{n-1}(t)) \neq 0, \,
\forall t \in (t_{0}-\varepsilon, t_{0}+\varepsilon) \}. $$
We assume that, $x(t) = \pi(\gamma(t))$,
where $\gamma(t)=\psi^{H+f}_{t}(\vartheta)$, does not contain autointersections
for $t \in (t_{0} -\varepsilon, t_{0} + \varepsilon)$, that is, that
$H_{p}(\gamma(t))=d\pi X^{H+f}(\gamma(t)) \neq 0$. Then there exists a
system of tubular coordinates $\mathcal{V}$,  in a neighborhood of
$\pi(\gamma(t_{0}))$, $F: \mathcal{V} \to \mathbb{R}^{n}$, such that,\\
i) $F(x)=(t,z_{1},...,z_{n-1})$;\\
ii) $F(x(t))=(t,0,...,0)$.\\
Observe that, by construction, $d_{x(t)}F H_{p}(\gamma(t))=(1,0,...0)$.
Consider a bump function $ \sigma :M \to \mathbb{R}$, such that, $\supp(\sigma) \subset \mathcal{V}$,  $\sigma \mid_{\mathcal{V}_{0}} \equiv 1$,  with $x(t_{0}) \in \mathcal{V}_{0} \subset \mathcal{V}$.
Define the perturbation space
$\mathcal{F}_{t_{0}} \subset C^{\infty}(M ; \mathbb{R})$ as being
$$\mathcal{F}_{t_{0}}=\{h_{\alpha,\beta}(x)=\tilde{h}_{\alpha,\beta}(x) \cdot \sigma(x)
 \, \mid \,\alpha, \beta \in \mathcal{A}_{t_{0}}\} $$
where, $\tilde{h}_{\alpha,\beta}(x)=\langle \alpha(t) \delta_{t_{0}}(t) + \beta(t)
\dot{\delta}_{t_{0}}(t) \, , \,  z \rangle$, $F(x)=(t,z)$ and
$\delta_{t_{0}}$ is a smooth approximation of the delta of Dirac in the point $t=t_{0}$.
Given $h_{\alpha,\beta} \in \mathcal{F}_{t_{0}}$ we get
$d_{x}h_{\alpha,\beta}=d_{x} \tilde{h}_{\alpha,\beta} \cdot \sigma(x) +
\tilde{h}_{\alpha,\beta} \cdot d_{x} \sigma(x)$. 
On the other hand
$$d_{x} \tilde{h}_{\alpha,\beta}=(\langle \frac{d \,}{dt}(\alpha(t)
\delta_{t_{0}}(t) + \beta(t)
\dot{\delta}_{t_{0}}(t)),z \rangle, \alpha(t) \delta_{t_{0}}(t) + \beta(t)
\dot{\delta}_{t_{0}}(t)) d_{x}F $$
Evaluating $x(t) $ and using that
$h_{\alpha,\beta}(x(t))=0$  and $\sigma(x(t))=1$, we get
$$d_{x(t)}h_{\alpha,\beta}=(0, \alpha(t) \delta_{t_{0}}(t) + \beta(t)
\dot{\delta}_{t_{0}}(t)) d_{x(t)}F.$$
In particular,
$$d_{x(t)}h_{\alpha,\beta} H_{p}(\gamma(t))=(0,\langle \alpha(t) \delta_{t_{0}}(t) + \beta(t)
\dot{\delta}_{t_{0}}(t)) d_{x(t)}F H_{p}(\gamma(t))=0$$
for any, $h_{\alpha,\beta} \in \mathcal{F}_{t_{0}}$.

We claim that,\\
1) $\mathcal{B}(\mathcal{F}_{t_{0}}) \subset T(H+f)^{-1}(k)$;\\
2) $X^{H+f}(\vartheta) \not\in \mathcal{B}(\mathcal{F}_{t_{0}})$;\\
3) $\dim(\mathcal{B}(\mathcal{F}_{t_{0}}))=2n-2$;\\
4) In particular, $\mathcal{B}(\mathcal{F}_{t_{0}})$
generates a space complementary to $\langle X^{H+f}(\vartheta), $ $Y^{H+f}(\vartheta)
\rangle $.\\
In order to get (1) consider,
 $\alpha_{0} = d_{x(t)}h_{\alpha,0}=(0, \alpha(t) \delta_{t_{0}}(t))
d_{x(t)}F = \alpha_{1}  \delta_{t_{0}}(t) $ and 
 $\beta_{0} = d_{x(t)}h_{0,\beta}=(0, \beta(t) \dot{\delta}_{t_{0}}(t))
d_{x(t)}F = \beta_{1} \dot{\delta}_{t_{0}}(t) $,
then,
$$   \mathcal{B}(h_{\alpha})= d_{\vartheta}\psi_{T}^{H + f}
\int_{0}^{T} (d_{\vartheta}\psi_{t}^{H + f })^{-1}
\left [
\begin{matrix}
0 \\
\alpha_{0}
\end{matrix}
\right ]
dt$$
\text{ and }   $$\mathcal{B}(h_{\beta})= d_{\vartheta}\psi_{T}^{H + f}
\int_{0}^{T} (d_{\vartheta}\psi_{t}^{H + f })^{-1}
\left [
\begin{matrix}
0 \\
\beta_{0}
\end{matrix}
\right ]
dt.$$
Observe that,
\( \omega( \left [  \begin{matrix}   0 \\   \alpha_{0}  \end{matrix}  \right
], X^{H+f}(\gamma(t)))=  \alpha_{0} H_{p}(\gamma(t))=0, \)
\text{ and }\\
\(\omega( \left [\begin{matrix}   0 \\   \beta_{0}  \end{matrix}  \right
], X^{H+f}(\gamma(t)))=  \beta_{0} H_{p}(\gamma(t))=0,\)
therefore
$\left [  \begin{matrix} 0 \\   \alpha_{0}  \end{matrix}  \right
]$ and $\left [  \begin{matrix}   0 \\   \beta_{0}  \end{matrix}  \right
]$ are in $T(H+f)^{-1}(k)$.
Thus, $\mathcal{B}(\mathcal{F}_{t_{0}})\subset
T(H+f)^{-1}(k)$.

In order to get (2),   we will make $\delta_{t_{0}} \to \delta_{Dirac}$
and will write $\mathcal{B}(h_{\alpha})$ and $\mathcal{B}(h_{\beta})$ as 
$$ \displaystyle \mathcal{B}(h_{\alpha})= d_{\vartheta}\psi_{T}^{H + f}(d_{\vartheta}\psi_{t_{0}}^{H + f
})^{-1} \left [ \begin{matrix} 0 \\ \alpha_{1}(t_{0}) \end{matrix} \right].$$
Analogously,  
$$ \displaystyle \mathcal{B}(h_{\beta})= d_{\vartheta}\psi_{T}^{H + f}(d_{\vartheta}\psi_{t_{0}}^{H + f
})^{-1}
\left\{ J \mathcal{H}^{H+f}(t_{0}) \left[ \begin{matrix} 0 \\ \beta_{1}(t_{0})
\end{matrix} \right] - \left[ \begin{matrix} 0 \\ \dot{\beta}_{1}(t_{0})
\end{matrix} \right]
\right\}. $$

If we assume, by contradiction, that $X^{H+f}(\gamma(t))=\mathcal{B}(h_{\alpha})
+\mathcal{B}(h_{\beta})$ then
$$X^{H+f}(\gamma(t_{0}))= \left\{ \left [ \begin{matrix}
0 \\ \alpha_{1}(t_{0}) \end{matrix} \right] +
J \mathcal{H}^{H+f}(t_{0}) \left [ \begin{matrix} 0 \\ \beta_{1}(t_{0})
\end{matrix} \right ] - \left [ \begin{matrix} 0 \\ \dot{\beta}_{1}(t_{0})
\end{matrix} \right ] \right\}.$$
From this equality we have
$H_{p}(\gamma(t_{0}))= H_{pp}(\gamma(t_{0})) \beta_{1}(t_{0})$.
Since $H_{p}(\gamma(t_{0})) \neq 0$   we have $n-1$ choices, linearly independent, for $\beta_{1}(t_{0})$.  Indeed, $dF$ is an isomorphism and for all $\beta (t_{0}) \in  \mathbb{R}^{n-1}$ we have $ \beta_{1}(t_{0})H_{p}(\gamma(t_{0}))=(0,\beta (t_{0}))dFH_{p}(\gamma(t_{0}))=0$.  Thus,
 $0=\beta_{1}(t_{0})H_{p}(\gamma(t_{0}))=\beta_{1}(t_{0})
H_{pp}(\gamma(t_{0})) \beta_{1}(t_{0})$,  contradicting the superlinearity of $H$.
For (3) observe that, in (2) we got the limit representation
$$\mathcal{B}(h_{\alpha})
+\mathcal{B}(h_{\beta})=  d_{\vartheta}\psi_{T}^{H + f}
(d_{\vartheta}\psi_{t_{0}}^{H + f })^{-1} $$ $$ \left\{ \left [ \begin{matrix} 0 &
\quad H_{pp}(\gamma(t_{0})) \\
I_{n} & -H_{xp}(\gamma(t_{0})) \end{matrix} \right]
\left [ \begin{matrix} \alpha_{1}(t_{0}) \\ \beta_{1}(t_{0})  \end{matrix} \right ]-
\left [ \begin{matrix} 0 \\ \dot{\beta}_{1}(t_{0})  \end{matrix} \right ] \right\}.$$
From this equation we get  $\dim( \{ \mathcal{B}(h_{\alpha})
+\mathcal{B}(h_{\beta}) \})= \dim(\{\alpha_{1}(t_{0}),$  $ \beta_{1}
(t_{0})\}=2n-2$, since $\left [ \begin{matrix} 0 & \quad H_{pp}(\gamma(t_{0})) \\
I_{n} & -H_{xp}(\gamma(t_{0})) \end{matrix} \right]$ is an isomorphism.

Finally, we observe that, the claim (1) is true independently of the
approximation $\delta_{t_{0}}$ of the delta of Dirac in the point  $t=t_{0}$.
Moreover the claims (2) and (3) still true for $\delta_{t_{0}}$, close enough to the delta of Dirac.
\end{proof}

%###################################################################################

The next theorem allow  us to  make a local perturbation  of a
periodic orbit nondegenerate of order $\leq m$ in such way that it becomes
nondegenerate of order $\leq 2m$. The proof is just for dimension 2
and the $n$-dimensional case is still open. Almost all th parts of the argument
are true in the $n$-dimensional case, but we do not know how to show
the surjectivity of the representation in this case.

\begin{theorem}\label{PerturbacaoLocaldaOrbita}\textnormal{(Local perturbation of periodic orbits)}
Let  $dim(M)=2$, $\mathbb{H}:~T^{*}M \to \mathbb{R}$ be a smooth, convex
and superlinear Hamiltonian and $\gamma= \{
\psi_{t}^{\mathbb{H}}(\vartheta_{0}) \mid 0 \leq t \leq T \}
\subseteq \mathbb{H}^{-1}(k)$, where $\mathbb{H}^{-1}(k)$ is a
regular energy level, $T$ is the minimal period of $\gamma$, and
$\gamma$ is isolated in this energy level, nondegenerate of order $
\leq m \in \mathbb{N}$. Then there exists a potential  $f_{0} \in
C^{\infty}(M,\mathbb{R})$ arbitrarily close to zero, with
$\supp(f_{0}) \subset \mathcal{U} \subset M$ such that, 
$\gamma$ is nondegenerate of order $\leq 2m$ to
$\mathbb{H}+f_{0}$. More over, $\mathcal{U}$ can be chosen arbitrarily small.
\end{theorem}

\begin{proof}
Choose $t_{0} \in (0,T)$ and $E(0)=\{e_{1}(0), e_{2}(0),
e_{1}^{*}(0), e_{2}^{*}(0) \}$ a symplectic frame in $\gamma(t_{0})$
with $e_{1}(0)=X^{\mathbb{H}}(\gamma(t_{0}))$. Consider
$$E(t)=\{e_{1}(t), e_{2}(t), e_{1}^{*}(t), e_{2}^{*}(t) \}$$, where
$\xi(t)= d_{\gamma(t_{0})} \psi_{-t}^{\mathbb{H}} \xi, \quad \forall
\xi \in E(0)$, for $t \in (0,r)$ with $r >0$ arbitrarily small.

Then we can decompose the matrix of the differential of the flow, in
the base $E(0)$, $[d_{\gamma(t_{0})}
\psi_{T}^{\mathbb{H}}]_{E(0)}^{E(0)} \in Sp(2)$, as
$ [d_{\gamma(t_{0})} \psi_{T}^{\mathbb{H}}]_{E(0)}^{E(0)} =
[d_{\gamma(t_{0}-r)} \psi_{r}^{\mathbb{H}}]_{E(0)}^{E(r)} \cdot
[d_{\gamma(t_{0})} \psi_{T-r}^{\mathbb{H}}]_{E(r)}^{E(0)}.$ 
By construction we have that $[d_{\gamma(t_{0}-r)}
\psi_{r}^{\mathbb{H}}]_{E(0)}^{E(r)}=I_{4}$, therefore 
$$  \quad \quad  \quad   \quad  \quad \quad  \quad \quad [d_{\gamma(t_{0})} \psi_{T}^{\mathbb{H}}]_{E(0)}^{E(0)} =
[d_{\gamma(t_{0})} \psi_{T-r}^{\mathbb{H}}]_{E(r)}^{E(0)}.  \quad\quad \quad   \quad \quad   \quad \quad   \quad 
\quad (1)$$
Consider $U$ an arbitrarily small neighborhood of
$\gamma(t_{0})$ in $T^{*}M$ and $r$ small enough, in  such way that,
$\hat{\gamma}= \{ \psi_{t}^{\mathbb{H}}(\vartheta_{0}) \mid t \in
(t_{0}-r,t_{0}) \} \subseteq \mathbb{H}^{-1}(k) \cap U$. Fix
$t_{1} \in (t_{0}-r,t_{0})$ and $V$ a neighborhood of
$\gamma(t_{1})$ in $T^{*}M$, small enough, in  such way that, $V \subset U$
and that $\gamma(t_{0}),\gamma(t_{0}-r) \not\in \overline{V}$.
Suppose  that we have $\tilde{\mathbb{H}}: T^{*}M \to \mathbb{R}$ a
smooth Hamiltonian re\-pre\-sen\-ting a smooth perturbation of
$\mathbb{H}$, such that, $\supp(\tilde{\mathbb{H}}-\mathbb{H})
\subset V$ and that $jet_{1}(\tilde{\mathbb{H}})\mid_{\gamma(t)}=
jet_{1}(\mathbb{H})\mid_{\gamma(t)}$, and
$ [d_{\gamma(t_{0})} \psi_{T}^{\tilde{\mathbb{H}}}]_{E(0)}^{E(0)} =
[d_{\gamma(t_{0}-r)} \psi_{r}^{\tilde{\mathbb{H}}}]_{E(0)}^{E(r)}
\cdot [d_{\gamma(t_{0})}
\psi_{T-r}^{\tilde{\mathbb{H}}}]_{E(r)}^{E(0)}$.
Since $\supp(\tilde{\mathbb{H}}-\mathbb{H}) \subset V$, we have 
$[d_{\gamma(t_{0})} \psi_{T-r}^{\tilde{\mathbb{H}}}]_{E(r)}^{E(0)}=
[d_{\gamma(t_{0})} \psi_{T-r}^{\mathbb{H}}]_{E(r)}^{E(0)}$. By (1),
$ [d_{\gamma(t_{0})}
\psi_{T-r}^{\mathbb{H}}]_{E(r)}^{E(0)}= [d_{\gamma(t_{0})}
\psi_{T}^{\mathbb{H}}]_{E(0)}^{E(0)}$, so
$$\quad \quad  \quad   \quad \quad \quad[d_{\gamma(t_{0})} \psi_{T}^{\tilde{\mathbb{H}}}]_{E(0)}^{E(0)} =
[d_{\gamma(t_{0}-r)} \psi_{r}^{\tilde{\mathbb{H}}}]_{E(0)}^{E(r)}
\cdot [d_{\gamma(t_{0})} \psi_{T}^{\mathbb{H}}]_{E(0)}^{E(0)} \quad  \quad \quad\quad
\quad (2)$$
From the construction of the perturbation described above we have
that $[d_{\gamma(t_{0})}
\psi_{T}^{\tilde{\mathbb{H}}}]_{E(0)}^{E(0)}$ has the expression
$$
[d_{\gamma(t_{0})}
\psi_{T}^{\tilde{\mathbb{H}}}]_{E(0)}^{E(0)} =
\left [
\begin{matrix}
1   & \alpha & \sigma & \beta  \\
0   & A & \hat{\alpha} & B  \\
0   & 0 & 1 & 0  \\
0   & C & \hat{\beta} & D  \\
\end{matrix}
\right ]
\, \in Sp(2),
$$
because, the energy level in $\gamma(t_{0})$  and
$\gamma(t_{0}-r)$ is the same for $\tilde{\mathbb{H}}$ and
$\mathbb{H}$. Thus it is invariant by the action of the flow of
both Hamiltonians. Let $\hat{Sp(2)}$ be the following subgroup of $Sp(2)$,
$$
\hat{Sp(2)} =
\left \{
\left [
\begin{matrix}
1   & \alpha & \sigma & \beta  \\
0   & A & \hat{\alpha} & B  \\
0   & 0 & 1 & 0  \\
0   & C & \hat{\beta} & D
\end{matrix}
\right ]
\in SL(4)
\,
\left |
\,
\left [
\begin{matrix}
\hat{\alpha} \\
\hat{\beta}
\end{matrix}
\right ]
\right.
=
\left [
\begin{matrix}
A &  B  \\
C &  D
\end{matrix}
\right ]
J
\left [
\begin{matrix}
\alpha & \beta
\end{matrix}
\right ]^{*}
\right \}
$$
and consider the projection $\pi : \hat{Sp(2)} \to Sp(1)$ given by
$$
\pi \left(\left [
\begin{matrix}
1   & \alpha & \sigma & \beta  \\
0   & A & \hat{\alpha} & B  \\
0   & 0 & 1 & 0  \\
0   & C & \hat{\beta} & D
\end{matrix}
\right ]
\right )=
\left [
\begin{matrix}
A &  B  \\
C &  D
\end{matrix}
\right ],
$$
which  is a homomorphism of Lie groups.
Observe that $$[d_{\gamma(t_{0}-r)}
\psi_{r}^{\tilde{\mathbb{H}}}]_{E(0)}^{E(r)}, [d_{\gamma(t_{0})}
\psi_{T}^{\mathbb{H}}]_{E(0)}^{E(0)} \in \hat{Sp(2)}$$ and
$\det([d_{\gamma(t_{0})} \psi_{T}^{\tilde{\mathbb{H}}}]_{E(0)}^{E(0)} -
\lambda I_{4})=(\lambda -1)^{2} \det(\pi([d_{\gamma(t_{0})}
\psi_{T}^{\tilde{\mathbb{H}}}]_{E(0)}^{E(0)}) - \lambda I_{2}).
$ Thus $\gamma$ will be a nondegenerate
orbit of order $\leq 2m$, to the perturbed Hamiltonian, if
$\pi([d_{\gamma(t_{0})}
\psi_{T}^{\tilde{\mathbb{H}}}]_{E(0)}^{E(0)})$ does not have roots
of the unity of order~$\leq~2m$ as eigenvalues. Since, the
symplectic matrices that are $2m$-elementary~\footnote{A symplectic
matrix is N-elementary if its principal eigenvalues (the
eigenvalues $\lambda$ such that $\| \lambda \| < 1$ or $Re(\lambda)
\geq 0$) are multiplicatively independent over the integer, that is,
if $\Pi \lambda_{i}^{p_{i}}=1$, where $\sum p_{i} =N$ then $p_{i}=0,
\, \forall i $.} (in particular, does not have roots of the unity of
order~$\leq~2m$ as eigenvalues), forms an open and dense subset of
$Sp(1)$,  we must to show that, for a choice of the perturbation
space , the correspondence
$\tilde{\mathbb{H}} \to \pi([d_{\gamma(t_{0})}
\psi_{T}^{\tilde{\mathbb{H}}}]_{E(0)}^{E(0)})$ applied to a
neighborhood of $\mathbb{H}$, generate  an  open neighborhood of
$\pi([d_{\gamma(t_{0})} \psi_{T}^{\tilde{\mathbb{H}}}]_{E(0)}^{E(0)})$ in $Sp(1)$.
Using the homomorphism property 
$$ \pi( [d_{\gamma(t_{0})} \psi_{T}^{\tilde{\mathbb{H}}}]_{E(0)}^{E(0)}) =
\pi( [d_{\gamma(t_{0}-r)}
\psi_{r}^{\tilde{\mathbb{H}}}]_{E(0)}^{E(r)}) \cdot \pi(
[d_{\gamma(t_{0})} \psi_{T}^{\mathbb{H}}]_{E(0)}^{E(0)} ). $$
We define $\mathbb{X}_{0}=\pi( [d_{\gamma(t_{0})}
\psi_{T}^{\mathbb{H}}]_{E(0)}^{E(0)} )$ and
$\hat{S}(\tilde{\mathbb{H}})=\pi( [d_{\gamma(t_{0}-r)}
\psi_{r}^{\tilde{\mathbb{H}}}]_{E(0)}^{E(r)})$.
Since the translation $\mathbb{X} \to \mathbb{X} \cdot
\mathbb{X}_{0}$ is an isomorphism of the of the Lie group $Sp(1)$,
we need to show that the map $\tilde{\mathbb{H}} \to
\hat{S}(\tilde{\mathbb{H}})$ applied to a neighborhood of
$\mathbb{H}$ generates an open neighborhood of $I_{2}$ in $Sp(1)$.
Inorder to construct the perturbation space we will consider
$\mathcal{N} \subset \mathbb{H}^{-1}(k)$ a local Lagrangian
submanifold in $\gamma(t_{0})$. We can reduce, if necessary, the
size of the neighborhood $U$ of $\gamma(t_{0})$ chosen previously
in such way that $U$ admits the parameterization
$(x=(x_{1},x_{2}),p=(p_{1},p_{2})):U \to \mathbb{R}^{2+2}$ as in
\cite{GonzaloPaternGenGeodFlowPositEntropy}, Lemma A3, that is,\\
a) $\mathcal{N} \cap U = \{(x,0) \}$; \\
b) $\omega=dx \wedge dp$;\\
c) $X^{\mathbb{H}}|_{\mathcal{N} \cap U}= 1 \frac{\partial}{\partial
x_{1}}$.\\
In these coordinates we can see that $\hat{\gamma}=\{ (t,0,0,0) \mid
t \in (t_{0}-r,t_{0}) \}$. Consider, the perturbation space,
$$\hat{\mathcal{F}}=\{ f:T^{*}M \to \mathbb{R} \mid \supp(f) \subset
\hat{W} \subset W \}$$ where $W = \mathcal{N} \cap V $ and $\hat{W}$
is a compact set contained in $W$ that contains $\gamma(t_{1})$ in its
interior.  Observe that, $\hat{\mathcal{F}}$ can be identified with
$C^{\infty}(\hat{W}, \mathbb{R})$, therefore we can think
$\hat{\mathcal{F}}$ as a vectorial space. Consider the
following finite dimensional subspace $\mathcal{F} \subset
\hat{\mathcal{F}}$
$$\mathcal{F}=\{ f \mid f(x,p)=\eta(x)(a \delta_{t_{1}}(x_{1}) + b
\delta_{t_{1}}'(x_{1}) + c
\delta_{t_{1}}''(x_{1}))\frac{1}{2}x_{2}^{2}, \, a,b,c \in
\mathbb{R} \}$$ where $\eta$ is a fix function with
$\supp(\eta) \subset \hat{W}$ and $\eta \equiv 1$ in some
neighborhood of $\gamma(t_{1})$, in $\mathcal{N}$. Moreover,
$\delta_{t_{1}}$ is a smooth approximation of the delta of Dirac in
the point $t_{1}$.  Now we are able to define the differentiable map
$S: \mathcal{F} \to Sp(1)$ given by
$S(f)=\hat{S}(\mathbb{H}+f)=\pi( [d_{\gamma(t_{0}-r)}
\psi_{r}^{\mathbb{H} +f }]_{E(0)}^{E(r)}) $. Observe that
$dim(\mathcal{F})=3=dim(Sp(1))$ and $S(0)=\hat{S}(\mathbb{H}+0)=
\pi( [d_{\gamma(t_{0}-r)} \psi_{r}^{\mathbb{H} +0 }]_{E(0)}^{E(r)})=
I_{2}$. Thus we must to show that,
$$d_{0}\mathcal{F} :T_{0} \mathcal{F} \cong \mathcal{F} \to
T_{Id_{2 \times 2}}Sp(1) \cong sp(1)$$ is surjective.
Given $h \in \mathcal{F}$ we have
$
d_{0}\mathcal{F}(h)=
\pi( \frac{d}{dl}  [d_{\gamma(t_{0}-r)} \psi_{r}^{\mathbb{H} +lh
}]_{E(0)}^{E(r)} |_{l=0} ).
$ 
Consider $\xi \in T_{\gamma(t_{0}-r)}T^{*}M$, where $t \in (0, r)$
 and define $$\xi(t,l)=d_{\gamma(t_{0}-r)}
\psi_{t}^{\mathbb{H} +lh } \xi.$$ For a fix $l$ we define a field
through $\gamma$ that verifies the equation
$$
\left \{
\begin{array}{clcr}
\dot{\xi}(t,l) = J Hess(\mathbb{H}+lh)(\gamma(t)) \xi(t,l) \\
\xi(0,l)=\xi. \quad \quad \quad \quad \; \quad \quad \quad \quad \quad \quad \; ^{ }
\end{array}
\right.
$$
Taking the derivative of the equation above with respect to $l$ and
using the commutativity of the derivatives we get
$$\frac{d}{dt}( \frac{d}{dl} \xi(t,l)|_{l=0}) = J Hess(h)
\xi(t,l)|_{l=0} + J Hess(\mathbb{H}) \frac{d}{dl} \xi(t,l)|_{l=0}$$
Denote  $\mathcal{H}= Hess(\mathbb{H})$, $\xi(t)=\xi(t,l)|_{l=0}$ e
$\mathbb{Y}(t)= \frac{d}{dl} \xi(t,l)|_{l=0}$, then 
$$
\left \{
\begin{array}{clcr}
\dot{\mathbb{Y}}(t) =J \mathcal{H} \mathbb{Y}(t)+ J Hess(h) \xi(t)\\
\mathbb{Y}(0)=0. \quad \quad \quad \quad \quad \quad \quad \quad \quad \quad \; ^{ }
\end{array}
\right.
$$
Applying the method of variation of constants and using
$$
\left \{
\begin{array}{clcr}
\dot{\xi}(t) = J \mathcal{H}(\gamma(t)) \xi(t) \\
\xi(0)=\xi, \quad \quad \quad \quad \quad \; ^{ }
\end{array}
\right.
$$
we get
$$\mathbb{Y}(t)=d_{\gamma(t_{0}-r)} \psi_{r}^{\mathbb{H}}
\int_{0}^{r} d_{\gamma(t_{0}-r)}\psi_{-t}^{\mathbb{H}}
J Hess(h) d_{\gamma(t_{0}-r)}\psi_{t}^{\mathbb{H}} \xi dt.$$
Remember that $\mathbb{Y}(r)= \frac{d}{dl} \xi(r,l)|_{l=0}=
\frac{d}{dl} d_{\gamma(t_{0}-r)} \psi_{r}^{\mathbb{H} +lh }(\xi)
|_{l=0}$, so
$$ \frac{d}{dl} d_{\gamma(t_{0}-r)} \psi_{r}^{\mathbb{H} +lh }
|_{l=0} =d_{\gamma(t_{0}-r)} \psi_{r}^{\mathbb{H}}
\int_{0}^{r} d_{\gamma(t_{0}-r)}\psi_{-t}^{\mathbb{H}}
J Hess(h) d_{\gamma(t_{0}-r)}\psi_{t}^{\mathbb{H}} dt.$$
From this calculation
$$\;  d_{0}\mathcal{F}(h)=\pi \left( [d_{\gamma(t_{0}-r)} \psi_{r}^{\mathbb{H}}
\int_{0}^{r} d_{\gamma(t_{0}-r)}\psi_{-t}^{\mathbb{H}} J Hess(h)
d_{\gamma(t_{0}-r)}\psi_{t}^{\mathbb{H}} dt]_{E(0)}^{E(r)}  \right ). \; 
  (3)$$
In order to obtain the expression (3) we need to calculate $J Hess(h)$.
All the integrals will be calculated  with the delta of
Dirac and not with the approximations, however the same conclusions
are true for an approximation, good enough.
Consider $\tilde{h}(x)=(a \delta_{t_{1}}(x_{1}) + b
\delta_{t_{1}}'(x_{1}) + c
\delta_{t_{1}}''(x_{1}))\frac{1}{2}x_{2}^{2}$ and
$h(x)=\eta(x) \tilde{h}(x)$ then
$d h= \eta d \tilde{h} + \tilde{h} d \eta$ and
$d^{2}h= \eta d^{2}\tilde{h} +d \tilde{h}^{*}  d \eta + d \eta^{*} d \tilde{h} + \tilde{h} d^{2}
\eta$. As, $Hess(h)(\gamma)= d_{\gamma}^{2}h$ and
$jet_{1}(h)|_{\gamma}=0$ we have that $Hess(h)(\gamma)=\eta(\gamma)
d_{\gamma}^{2}\tilde{h}$. On the other hand,
$(d_{\gamma}^{2}\tilde{h})_{ij}= a \delta_{t_{1}}(t_{0}-r +t) + b \delta_{t_{1}}'(t_{0}-r +t)
+ c \delta_{t_{1}}''(t_{0}-r +t)$ if $ij=22$ and, and equal to 0 otherwise.
Taking the $x_{1}$-support of $\delta_{t_{1}}$ small enough, we can
assume that
$
J Hess(h)(\gamma)=
\hat{A} \delta_{t_{1}}(t_{0}-r +t) +
\hat{B} \delta_{t_{1}}'(t_{0}-r +t) +
\hat{C} \delta_{t_{1}}''(t_{0}-r +t)
$
where 
$(\hat{A})_{ij}= -a $ if $ij=42$, and equal to 0 otherwise,
$(\hat{B})_{ij}= -b $ if $ij=42$, and equal to 0 otherwise and
$(\hat{C})_{ij}= -c $ if $ij=42$, and equal to 0 otherwise.\\
Denote,
\begin{align*}
\hat{I}_{1}= d_{\gamma(t_{0}-r)} \psi_{r}^{\mathbb{H}} \int_{0}^{r}
d_{\gamma(t_{0}-r)}\psi_{-t}^{\mathbb{H}} \; \hat{A} \;
d_{\gamma(t_{0}-r)}\psi_{t}^{\mathbb{H}} \; \delta_{t_{1}}(t_{0}-r +t) dt, \\
\hat{I}_{2}= d_{\gamma(t_{0}-r)} \psi_{r}^{\mathbb{H}} \int_{0}^{r}
d_{\gamma(t_{0}-r)}\psi_{-t}^{\mathbb{H}} \; \hat{B} \;
d_{\gamma(t_{0}-r)}\psi_{t}^{\mathbb{H}} \; \delta_{t_{1}}'(t_{0}-r +t) dt, \\
\hat{I}_{3}= d_{\gamma(t_{0}-r)} \psi_{r}^{\mathbb{H}} \int_{0}^{r}
d_{\gamma(t_{0}-r)}\psi_{-t}^{\mathbb{H}} \; \hat{C} \;
d_{\gamma(t_{0}-r)}\psi_{t}^{\mathbb{H}} \; \delta_{t_{1}}''(t_{0}-r +t) dt. \\
\end{align*}
Thus\\
$\hat{I}_{1} = d_{\gamma(t_{0}-r)} \psi_{r}^{\mathbb{H}} \;
d_{\gamma(t_{0}-r)}\psi_{-(t_{1}-t_{0}+r)}^{\mathbb{H}} \; \hat{A}
\; d_{\gamma(t_{0}-r)}\psi_{(t_{1}-t_{0}+r)}^{\mathbb{H}}, $\\
$\hat{I}_{2} = - d_{\gamma(t_{0}-r)} \psi_{r}^{\mathbb{H}} \;
d_{\gamma(t_{0}-r)}\psi_{-(t_{1}-t_{0}+r)}^{\mathbb{H}}
 \; [\hat{B}, J \mathcal{H}] \;
d_{\gamma(t_{0}-r)}\psi_{(t_{1}-t_{0}+r)}^{\mathbb{H}},
$\\
and\\
$\hat{I}_{3}=d_{\gamma(t_{0}-r)} \psi_{r}^{\mathbb{H}} \;
d_{\gamma(t_{0}-r)}\psi_{-(t_{1}-t_{0}+r)}^{\mathbb{H}} \;(
[[\hat{C}, J \mathcal{H}], J \mathcal{H}] + [\hat{C}, J
\dot{\mathcal{H}}] ) \; \\
d_{\gamma(t_{0}-r)}\psi_{(t_{1}-t_{0}+r)}^{\mathbb{H}}$.\\
Define $\mathcal{Z}= \hat{A} - [\hat{B}, J \mathcal{H}] +
[\hat{C}, J \dot{\mathcal{H}}] + [[\hat{C}, J \mathcal{H}], J
\mathcal{H}]$. Then,

$$d_{0}\mathcal{F}(h)=\pi \left( [d_{\gamma(t_{0}-r)} \psi_{r}^{\mathbb{H}} \;
d_{\gamma(t_{0}-r)}\psi_{-(t_{1}-t_{0}+r)}^{\mathbb{H}} \;
\mathcal{Z} \;
d_{\gamma(t_{0}-r)}\psi_{(t_{1}-t_{0}+r)}^{\mathbb{H}}]_{E(0)}^{E(r)}
\right ).$$
Writing this matrix in the bases $E(0)$ and $E(r)$,  
in each point of the curve, we get,\\
$$[d_{\gamma(t_{0}-r)} \psi_{r}^{\mathbb{H}} \;
d_{\gamma(t_{0}-r)}\psi_{-(t_{1}-t_{0}+r)}^{\mathbb{H}} \;
\mathcal{Z} \;
d_{\gamma(t_{0}-r)}\psi_{(t_{1}-t_{0}+r)}^{\mathbb{H}}]_{E(0)}^{E(r)}=$$
$$
=[d_{\gamma(t_{0}-r)} \psi_{r}^{\mathbb{H}}]_{E(0)}^{E(r)} \;
[d_{\gamma(t_{0}-r)}\psi_{-(t_{1}-t_{0}+r)}^{\mathbb{H}}]_{E(r)}^{E(t_{0}-t_{1})}
\; [\mathcal{Z}]_{E(t_{0}-t_{1})}^{E(t_{0}-t_{1})}$$ \\ $$\;
[d_{\gamma(t_{0}-r)}\psi_{(t_{1}-t_{0}+r)}^{\mathbb{H}}]_{E(t_{0}-t_{1})}^{E(r)}.$$
Moreover $[d_{\gamma(t_{0}-r)} \psi_{r}^{\mathbb{H}}]_{E(0)}^{E(r)}=
I_{4}$ and there exists a symplectic conjugation $\mathbb{G} \in $ 
$Sp(2)$ between the base $E(t_{0}-t_{1})$ and the canonic symplectic
base,
$\left\lbrace  \frac{\partial}{\partial x_{1}}(\gamma(t_{1})),
\frac{\partial}{\partial x_{2}}(\gamma(t_{1})), \right. $ 
$\left. \frac{\partial}{\partial p_{1}}(\gamma(t_{1})),
\frac{\partial}{\partial p_{2}}(\gamma(t_{1})) \right\rbrace $ such
that, \\ $[\mathcal{Z}]_{E(t_{0}-t_{1})}^{E(t_{0}-t_{1})} = 
\mathbb{G}^{-1} \mathcal{Z} \mathbb{G}$. Thus
 $d_{0}\mathcal{F}(h)=$
$$\pi \left( \mathbb{G} [d_{\gamma(t_{0}-r)}
\psi_{(t_{1}-t_{0}+r)}^{\mathbb{H}}]_{E(r)}^{E(t_{0}-t_{1})} \right)
^{-1} \pi \left( \mathcal{Z}\right) \, \pi \left( \mathbb{G}
[d_{\gamma(t_{0}-r)}\psi_{(t_{1}-t_{0}+r)}^{\mathbb{H}}]_{E(r)}^{E(t_{0}-t_{1})}
\right)$$ that is, we need to show that $\pi(\mathcal{Z})$ is
surjective in $sp(1)$. A simple calculation shows that,
$
\pi(\mathcal{Z})=
\left [
\begin{matrix}
z_{11}   & z_{12} \\
z_{21}   & z_{22} \\
\end{matrix}
\right ]
$
where,\\
$z_{11}=-b \mathbb{H}_{p_{2}p_{2}} +2 c \mathbb{H}_{p_{2}p_{2}}
\mathbb{H}_{x_{2}p_{2}} + \dot{\mathbb{H}}_{p_{2}p_{2}}$\\
$z_{12}= 2c (\mathbb{H}_{p_{2}p_{2}})^{2}$\\
$z_{21}=-a + 2 b \mathbb{H}_{x_{2}p_{2}}+ 2 c
\mathbb{H}_{p_{1}p_{2}} \mathbb{H}_{x_{1}x_{2}} + 2 c
\mathbb{H}_{p_{2}p_{2}}\mathbb{H}_{x_{2}x_{2}} -2 c
\mathbb{H}_{x_{2}p_{1}} \mathbb{H}_{x_{1}p_{2}}\\ -4 c
(\mathbb{H}_{x_{2}p_{2}})^{2} -2 c \dot{\mathbb{H}}_{x_{2}p_{2}}$\\
$z_{22}=b \mathbb{H}_{p_{2}p_{2}} -2 c \mathbb{H}_{p_{2}p_{2}}
\mathbb{H}_{x_{2}p_{2}} - \dot{\mathbb{H}}_{p_{2}p_{2}}$\\

Remember that,
$
sp(1)=
\left \{
\left [
\begin{matrix}
B   & C \\
A   & -B \\
\end{matrix}
\right ] \; | A,B,C,D \in
\mathbb{R}, \;
\right \}
$
and $\mathbb{H}_{p_{2}p_{2}} \neq 0$, thus we have the
surjectivity. In order to conclude the proof we must to find a potential in $M$ adapted to this perturbation. Consider $f \in \mathcal{F}$ arbitrarily close
to zero such that $\pi
([d_{\gamma(t_{0})}\psi_{T}^{\mathbb{H}+f}]_{E(0)}^{E(0)})$ is
nondegenerate of order $\leq 2m$. Let us remember that the
$x$-support of $f$ is contained in $W$ which is an arbitrarily small
neighborhood of $\gamma(t_{1})$ in $\mathcal{N}$.
Consider $(\hat{x},\hat{p})$ the canonic symplectic coordinates in
$\gamma(t_{1})$, and $\hat{\pi} : T^{*}M \to \mathbb{R}$ given by
$\hat{\pi}(\hat{x},\hat{p})=\hat{x}$. As we are free  to
dislocate the point $t_{1}$ by a $\varepsilon$ arbitrarily small, we
can use the twist property of the vertical fiber bundle as in
Lemma~\ref{PropTwistDoFibradoVertical} to conclude that
$\hat{\pi}|_{\mathcal{N}}$ is a local diffeomorphism in
$\gamma(t_{1})$. Take a diffeomorphism $q:W~\subset~\mathcal{N}~\to~M$ given by
$q(x)=\hat{x}$, where $(x,0)\equiv (\hat{x},\hat{p})$ in
$\mathcal{N}$. Choose the potential
$$f_{0}(\hat{x})=
 \left \{
 \begin{matrix}
 f(q^{-1}(\hat{x})) \quad \quad  x \in \hat{\pi}(W)\\
 0  \quad \quad  \quad \quad \quad \quad  x \not\in \hat{\pi}(W),
 \end{matrix}
 \right.
$$
by construction, $\mathbb{H}(\hat{x},\hat{p}) + f_{0}(\hat{x})$ has
the desired property. The lemma is proven.
\end{proof}

%###################################################################################
\textbf{Conjecture:}\\
\textit{If $\dim(M)=n$ then $\pi(\mathcal{Z})$ is  surjective in $sp(n-1)$.}\\

The main obstruction to prove this conjecture is that if $\mathcal{Z}= \hat{A} - [\hat{B}, J \mathcal{H}] +
[\hat{C}, J \dot{\mathcal{H}}] + [[\hat{C}, J \mathcal{H}], J
\mathcal{H}]$ then we need to solve equations like $UX+XU=D$, in the space of simetric $n-1 \times n-1$. But it is well known (see \cite{GonDomtd}) that, in our case, the solving of this type of equations requires additional hypothesis on the eigenvalues of $\mathbb{H}_{p p }$, which are not generic in Ma\~n\'e's sense. On the other hand, our approach it is essentially the only way to construct perturbations by potentials, thus we hope that in the future, we will be capable to solve this equation in higer dimensions.

%###################################################################################

\begin{lemma}\label{a2a E denso Em aa}
Given $k \in \mathbb{R}$, $f_{0} \in \mathcal{R}(k)$,
$\mathcal{U}_{f_{0}} \subseteq \mathcal{R}$ a $C^{\infty}$ neighborhood  of
$f_{0}$, $\alpha=\alpha(\mathcal{U}_{f_{0}})>0$, as in the
Lemma~\ref{MinimoPerioLema}, $a \in \mathcal{R}$ such that $0 < a <
\infty$ and $\mathcal{G}_{k}^{a,a} \cap \mathcal{U}_{f_{0}} \neq
\varnothing$, we have that $\mathcal{G}_{k}^{a,2a} \cap \mathcal{U}_{f_{0}}$ is dense in
$\mathcal{G}_{k}^{a,a} \cap \mathcal{U}_{f_{0}}$.
\end{lemma}
\begin{proof}
Take $f \in \mathcal{G}_{k}^{a,a} \cap \mathcal{U}_{f_{0}}$ and  $\mathcal{U}$ an arbitrary  
neighborhood of $f$. We must show that $\mathcal{U} \cap
(\mathcal{G}_{k}^{a,2a} \cap \mathcal{U}_{f_{0}}) \neq  \varnothing$.
From the definition of $\mathcal{G}_{k}^{a,a}$ we have that all periodic orbits
of  $H+f$ in the level $k$ with  minimal period $\leq
a$ are nondegenerate of order $m \leq \frac{a}{T_{min}}$.
Take $\rho_{f}: T^{*}M \times (0,a) \times (-\varepsilon,\varepsilon)
\to T^{*}M \times T^{*}M \times \mathbb{R}$. From 
Corollary~\ref{PeriodNivelSaoNaoDegenTrans} we have that
$\rho_{f} \pitchfork \Delta_{0}$.
Moreover, by Lemma~\ref{PeriodPontoSaoNaoDegenTrans},~(i),  
$$\rho_{f}^{-1}(\Delta_{0})=\{(\vartheta, T, 0) \mid \vartheta
\in (H+f)^{-1}(k), \; T \in (0,a) ,
\psi_{T}^{H+f}(\vartheta)=\vartheta \}$$
Observe that $\rho_{f}^{-1}(\Delta_{0}) \subset (H+f)^{-1}(k) \times
[0,a] \times \{0\}$ which is a compact set. As $\Delta_{0}$ is closed we have that
$\rho_{f}^{-1}(\Delta_{0})$ is a submanifold of dimension 1, with a finite
number of connected components. Since each periodic orbit,
$\{\psi_{t}^{H+f}(\vartheta) \mid
t \in [0,T], \; (\vartheta, T, 0) \in \rho_{f}^{-1}(\Delta_{0})\}
\subset \rho_{f}^{-1}(\Delta_{0})$, is a connected component of
dimension 1,the number of periodic orbits for
$H+f$ in the level $k$ with minimal period $\leq a$, distinct, is finite.
Denote, $\{\psi_{t}^{H+f}(\vartheta_{i}) \mid
t \in [0, T_{i} = T_{min}(\vartheta_{i})], \; (\vartheta_{i}, T_{i}, 0) \in
\rho_{f}^{-1}(\Delta_{0}) \}$, for $i=1,..N$, the $N$ periodic orbits for
$H+f$ in the level $k$, with its respective minimal periods.
From Theorem~\ref{PerturbacaoLocaldaOrbita} we can find a sum
of $N$ potentials $f_{0}=f_{1}+...+f_{N}$
arbitrarily close to 0, such that, all orbits are nondegenerate of order $ \leq 2m$
for $(H+f)+f_{0}$. The claim is proven because $f + f_{0} \in \mathcal{U} \cap
(\mathcal{G}_{k}^{a,2a} \cap \mathcal{U}_{f_{0}})$.
\end{proof}

%###################################################################################

\begin{lemma}\label{transa2a}
With the same notation of  the Lemma~\ref{a2a E denso Em aa}, if
$\mathcal{G}_{k}^{a,a} \cap \mathcal{U}_{f_{0}} \neq \varnothing$,
then $ev_{\rho} : \mathcal{G}_{k}^{a,2a} \cap \mathcal{U}_{f_{0}}
\times T^{*}M \times (0,2a) \times (-\varepsilon,\varepsilon)
\to T^{*}M \times T^{*}M \times \mathbb{R}$ is transversal to $\Delta_{0}$.
\end{lemma}
\begin{proof}
Indeed, given $(f, \vartheta, T, S) \in
\mathcal{G}_{k}^{a,2a} \cap \mathcal{U}_{f_{0}}
\times T^{*}M \times (0,2a) \times (-\varepsilon,\varepsilon)$,
if $ev(f, \vartheta , T, S) $ $\not\in \Delta_{0}$,
is done. So we can assume that $ev(f, \vartheta, T, S)
\in \Delta_{0}$, that is, $\vartheta$ is a periodic orbit of
 $H+f$ in the level $k$ with minimal period,
$T_{min}=T_{min}(\vartheta)$ and $S=0$.

If $T=T_{min}$ then $ev \pitchfork _{(f, \vartheta, T, 0)}
\Delta_{0}$ by Lemma~\ref{evaluationtransverlema}, (ii). On the
other hand, if $T =m T_{min},\; m \geq 2$ we have that $  m \leq  2a/T_{min}$, that is,
$ T_{min} \leq  2a /m \leq a$ so $\vartheta$ is nondegenerate of order $m$,
because$f \in \mathcal{G}_{k}^{a,2a} \cap \mathcal{U}_{f_{0}}$.
Thus $ev \pitchfork _{(f, \vartheta, T, S)} \Delta_{0}$ by
Lemma~\ref{evaluationtransverlema}, (i).
\end{proof}

%###################################################################################

\begin{lemma}\label{3meiosinter2a}
With the same notation of the Lemma~\ref{a2a E denso Em aa}, if
$\mathcal{G}_{k}^{a,a} \cap \mathcal{U}_{f_{0}} \neq \varnothing$,
then $(\mathcal{G}_{k}^{3a/2,3a/2} \cap
\mathcal{U}_{f_{0}}) \cap (\mathcal{G}_{k}^{a,2a} \cap \mathcal{U}_{f_{0}})$
is dense in $\mathcal{G}_{k}^{a,2a} \cap \mathcal{U}_{f_{0}}$.
\end{lemma}

\begin{proof}
Consider $\mathcal{B}= \mathcal{G}_{k}^{a,2a} \cap \mathcal{U}_{f_{0}}$,
which is a submanifold of $C^{\infty}(M;\mathbb{R})$ because it is open.
From the Lemma~\ref{transa2a} we have that  $ev_{\rho} : \mathcal{G}_{k}^{a,2a}
\cap \mathcal{U}_{f_{0}} \times T^{*}M \times (0,2a) \times (-\varepsilon,\varepsilon)
\to T^{*}M \times T^{*}M \times \mathbb{R}$ is transversal to $\Delta_{0}$.
Then,  Theorem~\ref{abraham} implies that $\mathfrak{R}=\{ f \in
\mathcal{G}_{k}^{a,2a} \mid \rho_{f} \pitchfork \Delta_{0} \}$ is a generic
subset of $\mathcal{G}_{k}^{a,2a} \cap \mathcal{U}_{f_{0}}$. In
particular $\mathfrak{R}$ is dense in
$\mathcal{G}_{k}^{a,2a} \cap \mathcal{U}_{f_{0}}$.
We claim that $\mathfrak{R} \subset (\mathcal{G}_{k}^{3a/2,3a/2} \cap
\mathcal{U}_{f_{0}}) \cap (\mathcal{G}_{k}^{a,2a} \cap
\mathcal{U}_{f_{0}})$.
Indeed, take $f \in \mathfrak{R}$, from Corollary~\ref{PeriodNivelSaoNaoDegenTrans},
all periodic orbits of the flow defined by $H+f$ in the energy level $k$,
with minimal period $T_{min}$ are nondegenerate of order
$m \leq \frac{2a}{T_{min}}$ because $\rho_{f} \pitchfork \Delta_{0}$.
If we have a periodic orbit for $H+f$ in the level $k$,
with minimal period $T_{min} \leq 3a/2$ take
$m' \leq \frac{3a/2}{T_{min}} = \frac{2a}{4/3 T_{min}} \leq
\frac{2a}{T_{min}}$ then this orbit is nondegenerate of order $m'$
in particular $f \in \mathcal{G}_{k}^{3a/2,3a/2} \cap \mathcal{U}_{f_{0}}$.
Therefore $(\mathcal{G}_{k}^{3a/2,3a/2}\cap \mathcal{U}_{f_{0}}) \cap
(\mathcal{G}_{k}^{a,2a}\cap \mathcal{U}_{f_{0}})$ is dense in
$\mathcal{G}_{k}^{a,2a} \cap \mathcal{U}_{f_{0}}$.
\end{proof}

%###################################################################################

\begin{lemma}\label{3a23a2 denso Em aa}
With the same notation of the Lemma~\ref{a2a E denso Em aa}, if
$\mathcal{G}_{k}^{a,a} \cap \mathcal{U}_{f_{0}} \neq \varnothing$,
we have that $(\mathcal{G}_{k}^{3a/2,3a/2} \cap
\mathcal{U}_{f_{0}})$ is dense in
$\mathcal{G}_{k}^{a,a} \cap \mathcal{U}_{f_{0}}$.
\end{lemma}
\begin{proof}
From Lemma~\ref{3meiosinter2a} we have that
$(\mathcal{G}_{k}^{3a/2,3a/2} \cap \mathcal{U}_{f_{0}}) \cap
(\mathcal{G}_{k}^{a,2a} \cap \mathcal{U}_{f_{0}})$ is dense in
$\mathcal{G}_{k}^{a,2a} \cap \mathcal{U}_{f_{0}}$ and from
Lemma~\ref{a2a E denso Em aa},  $\mathcal{G}_{k}^{a,2a} \cap \mathcal{U}_{f_{0}}$
is dense in $\mathcal{G}_{k}^{a,a} \cap \mathcal{U}_{f_{0}}$ therefore
$\mathcal{G}_{k}^{3a/2,3a/2} \cap \mathcal{U}_{f_{0}}$ is dense in
$\mathcal{G}_{k}^{a,a} \cap
\mathcal{U}_{f_{0}}$.
\end{proof}

%###################################################################################

{\large \textbf{  Proof of the Lemma~\ref{ReduLocalDaPrimParte}:}} \label{ProvaLemaReducao}

%###################################################################################

\begin{proof}
Given $k \in \mathbb{R}$, $f_{0} \in \mathcal{R}(k)$,
$\mathcal{U}_{f_{0}} \subseteq \mathcal{R}$ a $C^{\infty}$ neighborhood 
of $f_{0}$, $\alpha=\alpha(\mathcal{U}_{f_{0}})>0$ as in
Lemma~\ref{MinimoPerioLema}. Take $c \in \mathbb{R}_{+}$, if  $c <
\alpha$ then $\mathcal{G}_{k}^{c,c} \cap
\mathcal{U}_{f_{0}}=\mathcal{U}_{f_{0}}$ by Lemma~\ref{MinimoPerioLema}.
 So we can assume that $c \in \mathbb{R}_{+}$, with  $c \geq \alpha > a >0$.

We claim that, $\mathcal{G}_{k}^{(\frac{3}{2})^{\ell}a,
(\frac{3}{2})^{\ell}a} \cap \mathcal{U}_{f_{0}}$ is dense in
$\mathcal{G}_{k}^{a,a} \cap \mathcal{U}_{f_{0}}$, $\forall \ell
\in \mathbb{N}$. The proof is by induction in $\ell$.

For $\ell =1$ observe that, $\mathcal{G}_{k}^{a,a} \cap
\mathcal{U}_{f_{0}}=\mathcal{U}_{f_{0}} \neq \varnothing $, because
$\alpha > a >0$. Therefore $\mathcal{G}_{k}^{\frac{3}{2} a,
\frac{3}{2} a} \cap \mathcal{U}_{f_{0}}$ is dense in
$\mathcal{G}_{k}^{a,a} \cap \mathcal{U}_{f_{0}}$ by
Lemma~\ref{3a23a2 denso Em aa}. 

Suppose that, $\mathcal{G}_{k}^{(\frac{3}{2})^{\ell}a,
(\frac{3}{2})^{\ell}a} \cap \mathcal{U}_{f_{0}}$ is dense in
$\mathcal{G}_{k}^{a,a} \cap \mathcal{U}_{f_{0}}$, with $\ell
\geq 1$. Then $\mathcal{G}_{k}^{(\frac{3}{2})^{\ell}a,
(\frac{3}{2})^{\ell}a} \cap \mathcal{U}_{f_{0}} \neq \varnothing
$, from the density, and taking $a'=(\frac{3}{2})^{\ell}a$, we will have
that $\mathcal{G}_{k}^{ \frac{3}{2} a' , \frac{3}{2} a'} \cap \mathcal{U}_{f_{0}}$ is dense in
$\mathcal{G}_{k}^{a',a'} \cap \mathcal{U}_{f_{0}}$ by Lemma~\ref{3a23a2 denso Em aa}.
So $\mathcal{G}_{k}^{(\frac{3}{2})^{\ell+1}a,
(\frac{3}{2})^{\ell+1}a} \cap \mathcal{U}_{f_{0}}$ is dense in
$\mathcal{G}_{k}^{a,a} \cap \mathcal{U}_{f_{0}}$  concluding
the proof of the claim.

Consider $\ell_{0}$, such that, $(\frac{3}{2})^{\ell_{0}}a >
c$. Then, $\mathcal{G}_{k}^{(\frac{3}{2})^{\ell_{0}}a,
(\frac{3}{2})^{\ell_{0}}a} \cap \mathcal{U}_{f_{0}} \subset
\mathcal{G}_{k}^{c,c} \cap \mathcal{U}_{f_{0}} \subset
\mathcal{G}_{k}^{a,a} \cap \mathcal{U}_{f_{0}}=\mathcal{U}_{f_{0}}$.
Since $\mathcal{G}_{k}^{(\frac{3}{2})^{\ell_{0}}a,
(\frac{3}{2})^{\ell_{0}}a} \cap \mathcal{U}_{f_{0}}$ is dense in
$\mathcal{G}_{k}^{a,a} \cap \mathcal{U}_{f_{0}}$, we conclude that
$\mathcal{G}_{k}^{c,c} \cap \mathcal{U}_{f_{0}}$ is dense in
$\mathcal{U}_{f_{0}}$. The lemma is proven.
\end{proof}
\vskip 10mm

\vspace{0.3cm} \textbf{{\Large Acknowledgements }} \vspace{0.3cm}

I wish to thank Dr. Gonzalo Contreras by the invitation to work in CIMAT, where part of this work was done. I thank to my thesis advisor, Dr. Artur O. Lopes for several important remarks and corrections, and  Patrick Z. A. by your help in translating this manuscript. This work is part
of my PhD thesis in Programa de  P\'os-Gradua\c{c}\~ao em Matem\'atica - UFRGS.

%###################################################################################

\vspace{0.6cm}
E-mail: oliveira.elismar@gmail.com


\begin{thebibliography}{999}

\bibitem{Ab} R. Abraham,
{\slshape ``Transversality in manifolds of mappings''}, Bull. Amer.
Math Soc. N$^{\underline{0}}$ 69,(1963), pp 470--474.

\bibitem{AbRb} R. Abraham and J. Robbin,
{\slshape ``Transversal mappings and flows''}, Benjamin, New York,
1967.

\bibitem{AbMardRat}  R. Abraham, J. E. Marsden, T. S. Ratiou,
{\slshape ``Manifolds, Tensor Analysis And Applications''}, Global
Analysis, Pure And Applied, Addison Wesley, London, 1983.

\bibitem{Anos}  D. V. Anosov,
{\slshape ``Generic properties of closed geodesics''}, IZV. Akad.
Nauk. SSSR Ser. Mat. \textbf{46}, N$^{\underline{0}}$4(1982), pp
675-709, 896.

\bibitem{KbMg} K. Burns and M. Gidea,
{\slshape ``Differential Geometry and Topology. With a view to
dynamical systems.''}, Studies in Advanced Mathematics. Chapman \&
Hall/CRC, Boca Raton, FL, 2005. x+389 pp.


\bibitem{GoReConvex} G. Contreras \& R. Iturriaga, {\slshape ``Convex
Hamiltonians without conjugate points''}, Ergod. Th.  \& Dynam. Sys,
(1999), 19 , 901-952.

\bibitem{GonDomtd} G. Contreras, {\slshape ``Geodesic flows with positive topological
entropy, twist maps and dominated splittings''}, preprint.

\bibitem{GoReMin} G. Contreras, R. Iturriaga,
{\slshape ``Global Minimizers of Autonomous Lagrangians''}, To
Appear. URL: http://www.cimat.mx/\%7egonzalo/.

\bibitem{GonzaloPaternGenGeodFlowPositEntropy} G. Conterras,
G. P. Paternain, {\slshape ``Genericity of geodesic flows with
po\-si\-ti\-ve topological entropy on $\textbf{S}^{2}$.''}, Journal
of Differential Geometry 61 (2002), N$^{\underline{0}}$ 1, 1-49.

\bibitem{PhilHart} P. Hartman,
{\slshape ``Ordinary differential equations''}, John Wiley \& Sons
Inc., New York, 1964.

\bibitem{Kg}  W. Klingenberg, {\slshape ``Lectures on Closed Geodesics''}, Grundlehren der
Mathematischen Wissenschaften, 230, Springer-Verlag, Berlin-New
York, 1978.

\bibitem{KgTk} Wilhelm Klingenberg and Floris Takens,
{\slshape ``Generic properties of geodesic flows''}, Math. Ann. 197
(1972), 323-334.

\bibitem{ManeGenProp} R. Ma\~n\'e, {\slshape ``Generic properties
and problems of minimizing measures of Lagrangian systems.''}
Nonlinearity 9, (1996), N$^{\underline{0}}$ 2, 273--310.

\bibitem{Jan} J. A. G. Miranda {\slshape ``Generic properties for magnetic
flows on surfaces.''},  2006 Nonlinearity 19, 1849-1874.

\bibitem{PaternainGeodFlows}  G. P. Paternain,
{\slshape ``Geodesic Flows''}, Progress in Mathemathics,
Birkh\"user, Boston, Vol.180, 1999.

\bibitem{Peix} M. M. Peixoto, {\slshape ``On an approximation theorem of Kupka and Smale''}, J.
Differential Equations 3, (1967), 214-227.

\bibitem{RobI} C. Robinson, {\slshape ``Generic properties of conservative systems''}, Amer. J.
Math. 92, (1970), 562-603.

\bibitem{RobII} C. Robinson, {\slshape ``Generic properties of conservative systems II''}, Amer. J.
Math. 92, (1970), 897-906.

\end{thebibliography}
\end{document}